\documentclass[preprint,12pt,authoryear]{elsarticle}

\usepackage{color,soul}
\usepackage{gensymb}
\usepackage{amsmath,amssymb,bm}
\usepackage{nomencl}
\usepackage{soul}
\usepackage{graphicx,subfig,float}  
\usepackage{geometry}
\usepackage{pdflscape,rotating}
\usepackage[normalem]{ulem}
\usepackage{lineno,hyperref,cleveref}
\usepackage{lipsum}
\usepackage{lscape}
\usepackage{caption}
\usepackage{subfloat}
\usepackage{booktabs,amsfonts,dcolumn}
\usepackage{subfig}
\usepackage{array}
\usepackage{multicol}
\usepackage{multirow}
\usepackage{xcolor}
\usepackage{pgfplots}
\usepackage{tikz}
\usepackage{hhline}
\usepackage{textcomp}
\usepackage{relsize}
\usepackage{mathtools}

\journal{Journal}

\makenomenclature

\begin{document}
\begin{frontmatter}

\title{Downscaling Using CDAnet Under Observational and Model Noises: The Rayleigh-B\'enard Convection Paradigm}

\author[KAUST]{Mohamad Abed El Rahman Hammoud}
\author[UCambridge,TexasAM]{Edriss S. Titi}
\author[KAUST]{Ibrahim Hoteit} 
\author[KAUST]{Omar Knio \corref{cor}}

\cortext[cor]{Corresponding author Email: omar.knio@kaust.edu.sa}

\address[KAUST]{King Abdullah University of Science and Technology, Thuwal, Saudi Arabia, 23955}
\address[UCambridge]{Department of Applied Mathematics and Theoretical Physics, University of Cambridge, Cambridge CB3 0WA, UK}
\address[TexasAM]{Department of Mathematics, Texas A \& M University, College Station, TX 77843, USA}

\begin{abstract}


Efficient downscaling of large ensembles of coarse-scale information is crucial in several applications, such as oceanic and atmospheric modeling. 
The \textit{determining form map} is a theoretical lifting function from the low-resolution solution trajectories of a dissipative dynamical system to their corresponding fine-scale counterparts.
Recently, \citet{HammoudJAMES2022} introduced CDAnet a physics-informed deep neural network as a surrogate of the \textit{determining form map} for efficient downscaling.
CDAnet was demonstrated to efficiently downscale noise-free coarse-scale data in a deterministic setting. 
Herein, the performance of well-trained CDAnet models is analyzed in a stochastic setting involving (i) observational noise, 
(ii) model noise, and (iii) a combination of observational and model noises. The analysis is performed employing the Rayleigh-B\'enard convection paradigm,
under three training conditions, namely, training with perfect, noisy, or downscaled data.
Furthermore, the effects of noises, Rayleigh number, and spatial and temporal resolutions of the input coarse-scale information on the downscaled fields are examined. 
The results suggest that the expected $\ell_2$-error of CDAnet behaves quadratically in terms of the standard deviations of the observational and model noises.
The results also suggest that CDAnet responds to uncertainties similar to the theorized and numerically-validated CDA behavior with an additional error overhead due to CDAnet being a surrogate model of the determining form map. 


\end{abstract}

\end{frontmatter}


\section{Introduction}
\label{sec:intro}

High-resolution representations of oceanic and atmospheric flows can be achieved using sophisticated simulations based on high-fidelity general circulation models~\citep{Chassignet2021}. 
Recent developments with respect to data assimilation techniques resulted in reliable forecasts that accurately represent the states of the ocean and atmosphere, of which only sparse measurements are available~\citep{Srinivas2019}.
However, high-resolution computational experiments are still a challenge, particularly because only uncertain coarse resolution measurements of the system state are generally available~\citep{Fox2019, Njuki2020}. 
To address such uncertainties, probabilistic forecasts are typically performed using ensembles of model realizations, however, this approach results in increased computational complexity that scales with the ensemble size~\citep{Chassignet2017}.

Several techniques have been developed to exploit low-resolution information about an underlying dissipative dynamical system to predict its fine-scale evolution, a process referred to as downscaling~\citep{Assefa2019}. 
These algorithms are classified as statistical or dynamical \citep{Hannah2022}. 
Statistical downscaling aims to develop statistical relations between the low- and high-resolution solution fields, which can then be used to map coarse-scale data to its higher resolution counterpart~\citep{Jha2015, LAFLAMME2016}.
These techniques have been deployed in multiple settings, for instance, to increase the resolution of precipitation fields \citep{Schmidli2006}
, and estimate the seasonal statistics with respect to the significant wave height \citep{Wang2010}. 
Dynamical downscaling constrains the evolution of the high-resolution dynamical system with coarse data \citep{Hannah2022}. 
It is commonly performed using nudging techniques in which a high-resolution numerical solution is nudged toward another that is only described by coarse-scale information either point by point (grid nudging), or through its discrete low-modes (spectral nudging)~\citep{Altaf2017}. 
Dynamical downscaling is widely used in various climate applications including atmospheric \citep{Caldwell2009} and oceanic \citep{Camus2014} global circulation models. 
Statistical and dynamical downscaling techniques have been compared in several studies, with the results showing comparable performances under particular conditions \citep{Kidson1998, Altaf2017, LEROUX2018, Srinivas2019}.

Continuous data assimilation (CDA) is a dynamical downscaling technique that introduces a nudging term in the governing equations to continuously drive coarse-scales of the downscaled solution toward the coarse-scale data \citep{Azouani2013}. 
When the coarse data are sufficiently resolved and the nudging parameters are suitably selected, the CDA-downscaled fine solution approaches the reference solution at an exponential rate in time \citep{Azouani2013}. 
Utilizing the CDA framework, one is able to establish the existence of a more general Lipschitz map from the low spatial resolution trajectories to their higher resolution counterparts, called the determining form map \citep{Foias2012, Foias_2014, Foias2017, Biswas2019}.
The determining form map is essentially a \emph{lifting} function from coarse spatial-scale solution trajectories of the dissipative system to their higher resolution counterparts. 
CDA has been theoretically investigated for different dynamical systems, including the $2$D Rayleigh-B\'enard (RB) convection \cite{Farhat2015b, Farhat2017}, 3D planetary geostrophic model \cite{FarhatLunasinTiti2016}, and surface quasi-geostrophic equations \cite{Jolly2019}.
CDA and the determining form map were further numerically analyzed for downscaling the 2D Navier-Stokes \cite{gesho_olson_titi_2016}, global circulation models \cite{Srinivas2019, Desamsetti2022}, and 2D Rayleigh-B\'enard convection \citep{Altaf2017, Farhat2018, HammoudCOMG2022, HammoudJAMES2022}.

Ensembles are of utmost importance, especially in application to atmospheric and oceanic flows, primarily due to their ability to effectively capture uncertainties arising from observations, initial conditions, and model parameterizations \citep{Oppenheimer2016}. 
Ensembles provide a robust framework for understanding inherent errors and limitations associated with measurement devices and techniques. 
By executing multiple simulations that incorporate diverse plausible observational uncertainties, ensembles enable analyzing a range of potential outcomes, thus enhancing the reliability of the predictions \citep{Wang2018}. 
Ensembles serve as an essential tool for quantifying uncertainties originating from model parameterizations, which inherently involve approximations of complex processes within the models \citep{Iskandarani2016}, and initial conditions wherein minor variations can yield substantial disparities in the subsequent evolution. 
Through an ensemble-based framework with predictions coming from distinct parameter settings, one can effectively investigate the ramifications of such uncertainties on model outputs
\citep{Thacker2015}. 
Ultimately, ensembles play a crucial role in informing decision-making and risk assessment by providing probabilistic forecasts and enhancing our understanding of complex atmospheric and oceanic processes \citep{Du2019}.

In practical settings, low-resolution information (e.g.\ global climate projections and reanalyses) are typically provided in terms of ensembles of possible scenarios or realizations, which reflect model and observational uncertainties \citep{Lahoz2014, Tandeo2020}. 
Generally, models are imperfect with errors arising from input uncertainties and poorly known parameterizations used to account for unresolved processes \citep{Rasp2018}.
Moreover, data are usually noisy because of measurement errors or adverse environmental conditions~\citep{Rasti2018}. 
To address the challenges arising from uncertainties, a common approach is to downscale large ensembles of coarse-scale data, however, this approach necessitates independently downscaling each coarse resolution realization on a fine computational mesh, leading to significant computational requirements~\citep{Li2009, Kunii1000Mem}. 
The associated computational burdens necessitate an efficient downscaling procedure whose performance in the face of observational and model uncertainties is well analyzed for its proper utilization in real-world scenarios.


Recently, CDAnet \citep{HammoudJAMES2022} has been introduced as a physics-informed deep neural network (PI-DNN) that is trained to perform computationally-efficient downscaling of dissipative fluid flows, and was tested within the framework outlined in the theory of the determining form map \citep{Foias2012, Foias2016, Foias2017}. 
CDAnet serves as a surrogate of the determining form map, from the coarse-scale solution trajectory to its higher resolution counterpart. 
Under identical downscaling conditions, and in a deterministic setting that does not involve data or model uncertainties, CDAnet performs similarly to CDA without the requirement to numerically solve the governing equations~\citep{HammoudJAMES2022}. 
CDAnet could then be directly leveraged for efficiently downscaling large ensembles of coarse data, potentially at a fraction of the cost of independent fine-resolution simulations.
Nevertheless, CDAnet was only analyzed in a noise-free deterministic setting necessitating a thorough analysis of its performance at downscaling large ensembles in an uncertain framework.

In the present study, the performance of CDAnet was investigated in cases involving observational, model and combined noises, where an ensemble of downscaled fields are generated and their errors with respect to the reference solution are examined. 
The study presents a thorough analysis comparing the performance of CDAnet in the presence of uncertainties during the inference phase for three different training scenarios.
We focus on the RB convection to analyze the performance of CDAnet in all three scenarios, which are outlined as follows.

In the first setting, CDAnet is trained with perfect noise-free data and without model perturbations, as in~\citep{HammoudJAMES2022}, thus resulting in the “best” surrogate of the determining form map. 
This enables a thorough numerical examination of CDAnet and allows a comparison with the theoretical guarantees; eg.~\citep{Bessaih2015}.
The surrogate is then used to individually downscale each coarse-scale ensemble member in the network's inference stage involving a coarse-scale solution trajectory that was not present during CDAnet's training. 
The resulting fine-scale ensemble is assessed against the reference trajectory by means of several error metrics. 
The results are also contrasted with predictions obtained using CDA as a benchmark describing the best achievable downscaling performance.

The aforementioned setting is then relaxed in the second and third scenarios, inspired by practical applications in which coarse-scale data may be noisy and model parameterizations being generally imperfect. 
In particular, the second scenario involves observational and model uncertainties during the training process. 
Within this scenario, the resulting CDAnet model is an imperfect representation of the determining form map. 
This CDAnet model is then employed in an network inference setting to downscale coarse-scale solution trajectories in the presence of observational and model noise, where the resulting ensemble of downscaled fields is systematically analyzed. 

Finally, the third scenario addresses the situation where the fine-scale reference solution is not available. 
Consequently, CDA-downscaled solution fields are used to supervise CDAnet's training.
However, to remain within a practical setting, CDA downscaling was performed with slightly inaccurate downscaling parameters to introduce some ``model'' errors in the downscaled solution. 
Observational and model uncertainties are then incorporated during the inference to analyze the performance of CDAnet and contrast it with the previous two scenarios.
In all cases, we analyze the errors in the downscaled ensemble, examine its statistical properties, and compare the predictions with those obtained under idealized conditions.


The remaining of the manuscript is organized as follows. 
Section \ref{sec:RBC} introduces the RB convection problem.
The CDA algorithm and its theoretical underpinnings are then discussed in Sections~\ref{sec:CDA_} and \ref{sec:DetForm}, respectively. 
Section \ref{sec:CDAnet} outlines the CDAnet machinery, and Section \ref{sec:noise} describes the experimental setting describing how noisy data and model parameters are generated, in addition to the metrics used to assess the performance of the trained PI-DNN. 
The results for the three downscaling scenarios are then discussed in Sections \ref{sec:results1}-\ref{sec:results3}. 
Finally, Section \ref{sec:conclusion} summarizes the main conclusions of this study.



\section{Rayleigh-B\'enard Convection}
\label{sec:RBC}

The considered dynamical problem consists of the motion of an incompressible Newtonian fluid within a $2$D periodic channel, with a length period of $L_x$ in the horizontal $x$-direction and height $L_y$. 
A temperature difference across the solid boundaries of the channel is assumed to be maintained between the warm bottom boundary and cooler top wall.
Buoyancy effects induced by this temperature gradient result in thermofluidic instability \citep{Pandey2018}. 
The onset of instability occurs at a critical Rayleigh number, ${\rm Ra}_c$, beyond which the RB convection is established. 
The RB flow is steady with well-defined cells for Rayleigh numbers, ${\rm Ra}$, slightly larger than ${\rm Ra}_c$, and the flow becomes chaotic with the increase in ${\rm Ra}$ ~\citep{Curry1984, Paul2007}. 

According to the scaling proposed by \citet{LEQUERE1991} and \citet{LEMAITRE2002}, the equations of motion employing the Boussinesq approximation can be expressed as follows:

\begin{equation}
    \nabla \cdot \mathbf{u} = 0,
    \label{eqn:incompr_Continuity}
\end{equation}
\begin{equation}
    \frac{\partial \mathbf{u}}{\partial t} +  \left( \mathbf{u} \cdot \nabla \right) \mathbf{u} + \nabla p = \frac{Pr}{\sqrt{Ra}} \nabla ^2 \mathbf{u} + Pr T \mathbf{e}_y,
    \label{eqn:momentum_true}
\end{equation}
\begin{equation}
    \frac{\partial T}{\partial t} + \left( \mathbf{u} \cdot \nabla \right) T = \frac{1}{\sqrt{Ra}} \nabla ^2 T + \mathbf{u} \cdot \mathbf{e}_2,
    \label{eqn:energy_true}
\end{equation}
where $T$, $p$, and $\mathbf{u} = \left( u, v \right)$ are the normalized temperature anomaly, pressure, and velocity fields, respectively;
the given physical parameters $Ra = \left(g \tilde{\alpha} \Delta \tilde{\Theta} \tilde{L}_y^3\right) \left(\tilde{\nu} \tilde{\kappa} \right)^{-1}$ is the Rayleigh number; 
$Pr = \nu \kappa ^{-1}$ is the Prandtl number; $\Delta \tilde{\Theta}$ is the absolute temperature difference between the plates; $g$ is the gravitational acceleration pointing downwards; and $\mathbf{e}_y$ is the unit vector in vertical direction.
We used $\tilde{\alpha}$, $\tilde{\nu}$, and $\tilde{\kappa}$ to denote the thermal expansion coefficient,  
kinematic viscosity, and thermal diffusivity, respectively, and tildes to denote dimensional quantities.

At the boundaries, the velocity and temperature are specified as follows:

\begin{equation}
    \mathbf{u}\left( t; x, 0\right) = \mathbf{u}\left( t; x, 1\right) = \mathbf{0}, 
\end{equation}
\begin{equation}
    T\left(t;x, 0 \right) = T\left(t;x, 1 \right) = 0.
\end{equation}
The equations are simulated subject to horizontally periodic boundary conditions in the 2D fundamental domain, $\Omega = \left[0, L_x/L_y \right] \times \left[0, 1 \right]$. For this purpose, the 
system is initialized with random velocity and temperature fields drawn from uniform distributions that are independent and spatially-uncorrelated.
Let $\mathcal{U}\left(a,b\right)$ denote a uniform distribution on the interval $(a,b)$,
then the initial conditions are specified as per $u\left( 0; x, y \right) \sim \mathcal{U}\left(-0.1,\ 0.1\right)$, 
$v\left( 0; x, y \right) \sim \mathcal{U}\left(-0.1,\ 0.1\right)$, and $T\left( 0; x, y \right) \sim \mathcal{U}\left(-0.1,\ 0.1\right)$ 
for all $\left(x, y\right) \in \Omega$.

Note that the above-mentioned initialization approach is motivated by the fact that the RB system is strongly sensitive to the initial conditions \citep{Chertovskih_2015, Yigit2015}. 
Thus, performing initializations with fields drawn from random distributions enables the construction of a rich dataset with diverse features, which can enhance the training of the CDAnet models.

\subsection{Numerical Setup}
\label{ssec:numSetup}

We follow the same numerical setup as in \citep{HammoudJAMES2022}, which is briefly described here.
A finite-difference methodology was used to solve the governing equations on a uniformly-spaced, staggered, Cartesian grid.
The domain was discretized into $n_x$ and $n_y$ points in the $x$ and $y$ directions, respectively, resulting in a computational grid with width $\Delta x = L_x/(n_xL_y)$ and height $\Delta y = L_y/n_y$. 
Furthermore, a second-order central differencing scheme is used to approximate spatial derivatives, and a third-order Adam-Bashforth scheme is selected for time integration. 
Finally, the continuity equation is satisfied using a pressure projection scheme~\citep{Chorin1968}, where the corresponding elliptic equation for pressure is efficiently solved using a fast Fourier transform algorithm.

The RB equations are solved on a 2D, periodic, 3:1 aspect ratio domain, i.e.,\ $L_x =  3$ and $L_y = 1$.
We used $n_x = 768$ points in the horizontal direction, $n_y = 256$ points in the vertical direction, and a fixed time step, $\Delta t = 5 \times 10^{-4}$.
Time integration was performed until a final simulation time, $t_f = 60$.
Three different Rayleigh numbers, $Ra=10^5$, $10^6$, and $10^7$ were considered with a fixed Prandtl number $Pr=0.7$ to observe different chaotic flow regimes~\citep{Lohse2010, plumley2016}.
Based on the aforementioned parameters, the grid Reynolds number was below 10 and the CFL number remained below $0.15$ throughout the simulation period. 
This ensured stability of the numerical solution, and that the simulated flow fields were suitably resolved.


\section{Continuous Data Assimilation}
\label{sec:CDA_}

CDA nudging is achieved by introducing into the governing equations a forcing term that is proportional to the difference between the spatial ``modes” interpolating the coarse data and those interpolating the high-resolution solution.
CDA guarantees the convergence of the downscaled solution at an exponential rate in time if the nudging term is suitably scaled and the coarse spatial-scale information has sufficient density.

Let $I_{h^o}$ be the spatial interpolation operator of an interpolant function $\phi$ specified on a uniform observation grid of size $h^o$. Then, 

\begin{equation}
    I_{h^o}\left( \phi \left( \mathbf{x} \right) \right) = \sum_{k = 1}^{N_{h^o}} \phi \left( \mathbf{x}_k \right) \chi_{Q_k} \left( \mathbf{x} \right),
\end{equation}
where $\mathbf{x}=(x, y)$, $Q_k$ are disjoint subsets such that diam($Q_k$) $\leq h^o$ for $k=1,\ ...,\ N_{h^o}$, $N_{h^o}$ is 
the number of observation points, $\bigcup_{j=1}^{N_{h^o}} Q_j = \Omega$, $\mathbf{x}_k \in Q_k$ are arbitrarily chosen points, and $\chi_S$ is the indicator function of set $S$. 
Let $\mathbf{w}=\left(\tilde u, \tilde v\right)$, $\Psi$, and $\varrho$ denote the CDA-downscaled velocity vector, temperature, and pressure, respectively. 
Then, the governing equations for the CDA-downscaled solution are governed by the system:

\begin{equation}
    \nabla \cdot \mathbf{w} = 0,
    \label{eqn:incompr_ContinuityCDA}
\end{equation}
\begin{equation}
    \frac{\partial \mathbf{w}}{\partial t} +  \left( \mathbf{w} \cdot \nabla \right) \mathbf{w} + \nabla \varrho = \frac{Pr}{\sqrt{Ra}} \nabla ^2 \mathbf{w} + Pr \Theta \mathbf{e}_y + \mu_{\mathbf{u}} \left( I_{h^o} \left( \mathbf{u}^o \right) - I_{h^o} \left( \mathbf{w} \right) \right),
    \label{eqn:momentum_CDA}
\end{equation}
\begin{equation}
    \frac{\partial \Psi}{\partial t} + \left( \mathbf{w} \cdot \nabla \right) \Psi = \frac{1}{\sqrt{Ra}} \nabla ^2 \Psi  + \mu_{T} \left( I_{h^o} \left( T^o \right) - I_{h^o} \left( \Psi \right) \right),
    \label{eqn:energy_CDA}
\end{equation}
where $\mu_{T}$ and $\mu_{\mathbf{u}}$ are strictly positive user-defined constants called the nudging parameters, and superscript “$o$" is used to represent an observed quantity. Here, we set $\mu_T = \mu_{\mathbf{u}}$ and tune them such that the downscaled solution rapidly converges to the reference solution.
Note that it was demonstrated theoretically \citep{FARHAT2015} and numerically \citep{Altaf2017} that downscaling using velocity observations alone is sufficient for recovering the RB convection solution. 
In this particular setting, $\mu_T$ would be set to zero and only velocity observations are used to nudge the system of equations. 
Furthermore, \citet{Altaf2017} showed that CDA's downscaled solution converges to the reference comparably regardless of the interpolant operator adopted.


\section{Determining Form Map}
\label{sec:DetForm}

For dissipative dynamical systems, the \textit{determining form} map is a Lipschitz continuous map from the trajectories of the coarse spatial-scales to the corresponding trajectories of high-resolution solution fields.
The determining form map establishes the theoretical grounds upon which the CDA algorithm was constructed.
The existence of this map was theoretically proven \citep{Foias2012, Foias_2014, Foias2017}, establishing the possibility of reducing the global attractor of the Navier-Stokes equations into an ordinary differential equation in the Banach space of the coarse spatial-scales trajectories \citep{Foias2017}. 
The determining form map was originally constructed for recovering the solution of the 2D Navier-Stokes equations by using a limited number of low-Fourier modes \citep{Foias2012}.
The framework was generalized by \citet{Foias_2014}, who demonstrated that the projection of the solution trajectories onto the coarse spatial-scales enables the construction of the determining form of any dissipative dynamical system.

\section{CDAnet}
\label{sec:CDAnet}

Although the determining form map provides an attractive framework for downscaling, the determining form map is generally an intractable function that does not admit a closed-form solution. 
Recently, however, \citet{HammoudJAMES2022} introduced CDAnet, an efficient, PI-DNN surrogate of the determining form map.
While it was thoroughly investigated in the noise-free deterministic setting, no study examines its performance in the uncertain setting.
The objective of this study is to examine the performance of CDAnet at downscaling large ensembles in the cases involving observational and model uncertainties.
All the experiments of this study were conducted using the same model architecture and training parameters as the one presented in \citep{HammoudJAMES2022}.
Consequently, the network architecture is briefly described, and readers are referred to the original paper for a fully detailed description on CDAnet's architecture and the training procedure. 
CDAnet is a two-stage PI-DNN comprising a feature extractor, followed by a task-performing multilayer perceptron (MLP).
The CDAnet architecture involves elements from the U-net \citep{unet2015}, and IMnet \citep{Chen_2019_CVPR}; however, it displays certain important differences, as discussed below.

The features extractor in CDAnet is a U-net convolution neural network \citep{unet2015} that takes the $3$D spatiotemporal low-resolution trajectory as input and extracts a $3$D features grid.
Moreover, CDAnet depends on inception blocks \citep{Szegedy2015, Szegedy2017} to magnify the receptive field and improve the quality of the extracted features.
The extracted features are represented by a grid comprising a features vector for each spatiotemporal grid coordinate point. 
These features describe the relationship of the solution variables over a spatiotemporal patch determined by the receptive field. 

The MLP then uses the features vector, which contains spatiotemporal information, and the corresponding spatiotemporal coordinates to predict the solution variables at that coordinate location.
In other words, the MLP predicts the solution variables on the higher-resolution grid on a point-by-point basis using spatiotemporal correlations between the input variables.
Similar to the work of \citet{Chen_2019_CVPR}, the features vector and spatiotemporal coordinates are concatenated with the hidden layers of the MLP, except for the final hidden layer.
CDAnet also includes the residual of the governing equations as an additional loss for training the PI-DNN to constrain the prediction of the PI-DNN, where this term is called the PDE loss and is scaled in the loss function using a positive coefficient $\lambda$.
Note that \citet{Chen_2019_CVPR} examined the proposed IMnet in the context of $3$D pointclouds, whereas CDAnet deals with $3$D spatiotemporal input and output data.

For the sake of completeness, we note that CDAnet share some commonalities with \textit{meshfreeflowNet} \citep{Max2020}, however, major differences are outlined below. 
The stark dissimilarity between CDAnet and \textit{meshfreeflowNet} is that CDAnet aims to identify a surrogate of the determining form map, whereas \textit{meshfreeflowNet} aims to compress a high-resolution solution trajectory.
In particular, \textit{meshfreeflowNet} is trained and evaluated on a common dataset composed of single or multiple solution trajectories, where no out-of-distribution samples are considered. 
This negates the need for the PDE loss, as seen by comparing the results of \citep{HammoudJAMES2022} and \citep{Max2020}.
Additionally, the convolution neural network of \textit{meshfreeflowNet} is built using ResNet blocks \citep{He_2016_CVPR}, which unlike CDAnet’s inception block, have a smaller receptive field \citep{Szegedy2017}.
Furthermore, CDAnet is trained using all the points composing the spatio-temporal grid in the training data instances, whereas \textit{meshfreeflowNet} relies on a number of randomly sampled points.

\section{Experimental Setup}
\label{sec:noise}

This section outlines the approaches used to generate noisy observations and perturbed model parameters.
Note that this study is primarily interested in the cases when observational and model noise are present during the inference stage, where a low-resolution trajectory that was not seen during training nor validation is downscaled in an uncertain setting. 
This involves three scenarios characterized by CDAnet's training, where training is noise-free, perturbed or comprised of imperfectly downscaled data. 

\subsection{Observational Noise}
\label{ssec:obsNoise}

We begin by simulating equations (\ref{eqn:incompr_Continuity})-(\ref{eqn:energy_true}) using the initial and boundary conditions discussed in Section \ref{sec:RBC}, resulting in a high-resolution reference solution. 
The high-resolution solution is then coarsened by subsampling the high-resolution fields at regular intervals using spatial and temporal downscaling factors $\mathcal{S}$ and $\mathcal{T}$, respectively. 
Here, $\mathcal{S}$ and $\mathcal{T}$ are positive integers that determine the coarsening factors used to generate the coarse-scale data relating the low-resolution observational grid resolution to the high-resolution solution grid. 
Compared to the original high-resolution reference solution, the resulting low-resolution data are $\mathcal{S}^2$ times coarser in space and $\mathcal{T}$ times coarser in time. 
Note that $\mathcal{S}$ and $\mathcal{T}$ are introduced for simplifying the analysis, but they should not be confused with the physical resolution of the observational grid on which the CDA's convergence relies upon.
In this study, we restrict the analysis to the cases of the spatial and temporal resolutions where the CDA downscaled solution converges to the reference in the noise-free deterministic setting. 


Noisy low-resolution trajectories are then obtained by adding independent and identically distributed (i.i.d.) white noise, randomly sampled from a Gaussian distribution with zero mean and standard deviation $\sigma_{obs}$.
Here, $\sigma_{obs}$ is fixed for all the solution variables, which amounts to perturbing the coarse-scale velocity components, temperature, and pressure with the same noise level.
While noisy observational data are considered during inference for all training scenarios, the study also considers a more practical case when observational noise is present during training. 

\subsection{Model Noise}
\label{ssec:modNoise}

Noisy model realizations are obtained by perturbing the weights of a neural network with random noise sampled from i.i.d. Gaussian distributions with zero mean. 
For each convolutional or dense layer, the noise level was defined in terms of a factor $\sigma_{mod}$ that scales the standard deviation of the unperturbed learnable parameters of the corresponding layer. 
In other words, perturbations are added to the weights of each convolution and dense layer by sampling $\mathcal{N}(0, \sigma_{mod}\times \sigma_{layer})$, where $\mathcal{N}$ denotes a Gaussian distribution, $\sigma_{mod}$ denotes the model noise level, and $\sigma_{layer}$ denotes the standard deviation of the weights associated with a given layer.
Note that $\sigma_{layer}$ ensures that the added noise scales appropriately in relation to the weights of a given layer.
In this setup, the weights of a convolution layer and a dense layer are assumed to be normally distributed with zero mean and a scaled standard deviation based on which the level of model errors can be described.
This strategy allows for analyzing a perturbed counterpart of the lifting function, where the noise added to the convolution layers impacts all spatial scales because the feature maps used to perform predictions would be perturbed.

Note that, this approach is in spirit of that used by \citet{Gagne2020}, which adds Gaussian noise to the outputs of each convolutional and dense layer of the network. 
The training of a CDAnet model with model noise is similar to the stochastic parameterization in a DNN \citep{Gagne2020}, where when the model parameters are perturbed, they inherit the perturbation's distribution.
The resulting CDAnet model becomes stochastic, where model parameters can be sampled from a distribution, leading to an ensemble of CDAnet models.
When model noise is considered during training, the optimization is affected by the perturbations, resulting in an uncertainties-aware CDAnet model.
While model noise is considered during inference for all training scenarios, the study also considers a more practical case when model noise is present during training.


\section{Numerical Results}

The results section is organized as three subsections corresponding to the different CDAnet training scenarios.
The first section considers perfect noise-free training conditions, the second examines the case of noisy training conditions and the third examines the case when imperfectly downscaled data are used for training.
Each of these sections is divided into three subsections characterizxing the noise involved during inference.
These subsections consist of CDAnet's inference when only observational noise, only model noise, or both observational and model noises are present.
CDAnet is then employed to downscale coarse-scale information about a test reference solution trajectory, which is not part of the training nor validation datasets. 
The resulting downscaled fields are compared to the noise-free reference by means of a number of error metrics to thoroughly analyze CDAnet's performance.
The error metrics include the mean absolute error (MAE), root mean squared error (RMSE), relative RMSE (RRMSE) and the expected $\ell_2$-error norm ($\Lambda$) and are outlined in \ref{ssec:evalMetrics}.

\subsection{Perfect Training Conditions}
\label{sec:results1}

In this section, CDAnet models are trained in the absence of observational and model errors resulting in the ``best'' surrogate of the determining form.
We analyze the impact of the presence of observational and model noise during the inference stage on CDAnet's downscaling performance. 
This idealized training setting enables a comparison between the theoretical guarantees \citep{Bessaih2015} and the present numerical results.
Furthermore, this settings allow comparing the performance of CDAnet to CDA in the face of uncertainties.

Note that multiple combinations of different $(\mathcal{S}, \mathcal{T})$ pairs and different Ra numbers were examined.
For brevity, we presented the results corresponding to $(\mathcal{S}, \mathcal{T}) = (4, 4)$ and $Ra=10^5$ only, while results for other parameters are included in the Supplementary Material.
Furthermore, for completeness, we note that the learning rate, batch size and all other training parameters are adapted from \citep{HammoudJAMES2022}.
For completeness, the learning rate, batch size and $\lambda$ are 0.2, 20 and 0.01, respectively.

\subsubsection{Effect of Observational Noise}
\label{res1_obs}

The impact of observational noise was investigated by downscaling 50 noisy realizations of a coarse-scale trajectory.
A wide range of observational noise levels was considered, such that $\sigma_{obs} \in  \left[0.001,\ 0.3\right]$.
Boxplots of the error metrics for the downscaled solution variables were first analyzed. 
The boxplots indicate that the RMSE and RRMSE of the downscaled solutions increase with increasing noise levels.
Furthermore, the results indicate flat boxplots meaning that the errors corresponding to different downscaled realizations are almost identical collapsing onto a single value.
This indicates that CDAnet behaves similarly to CDA in the case when only observational noise is present \citep{HammoudCOMG2022}. 

\begin{figure}
    \centering
    \includegraphics[width=\linewidth]{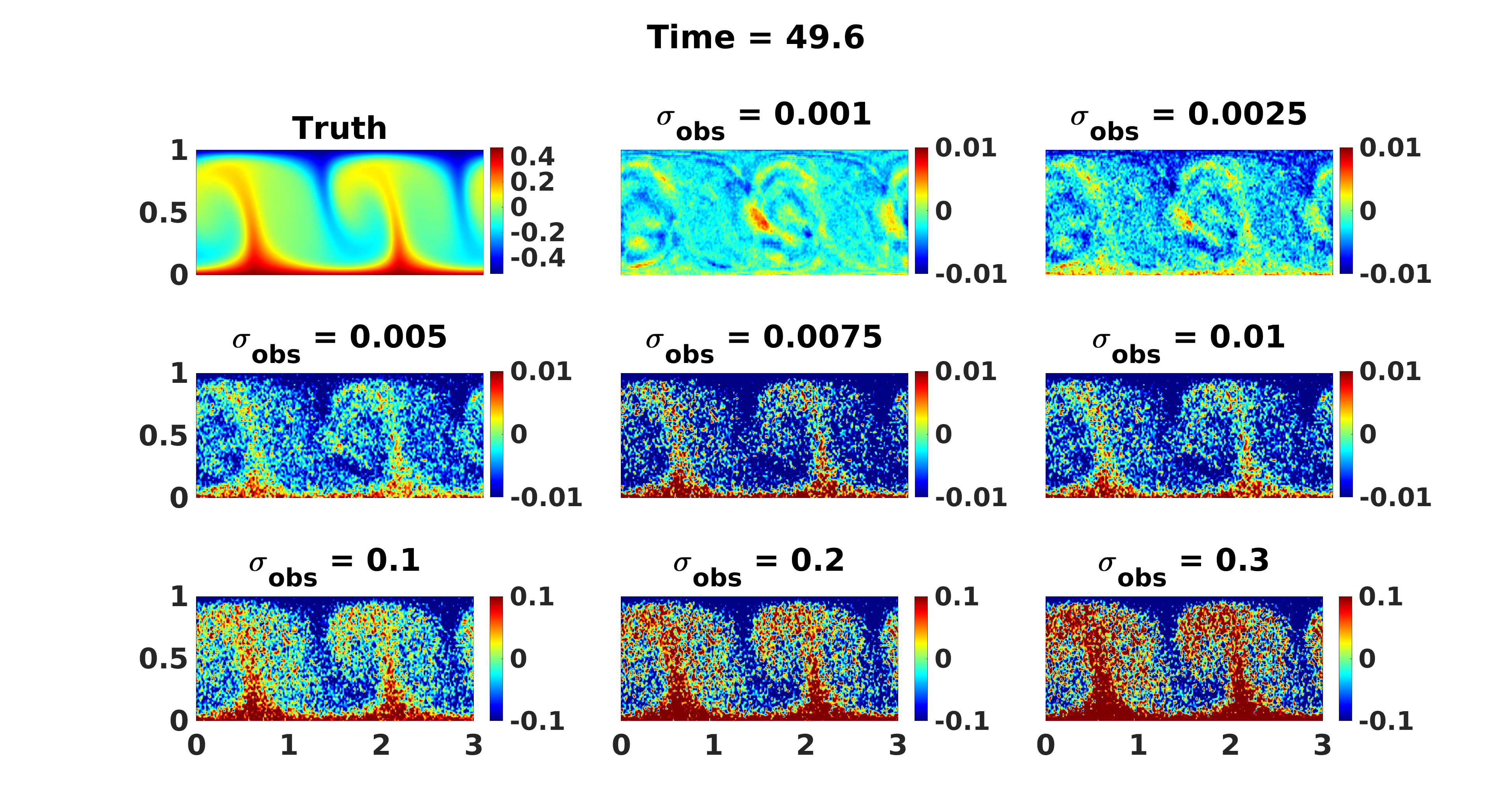}
    \caption{Snapshots of the temperature field corresponding to the reference solution, and the difference between a randomly selected downscaled realization and the reference solution for different values of $\sigma_{obs}$. The CDAnet model used to downscale the noisy coarse-scale trajectories was trained under noise-free condition for $\mathcal{S} = 4$, $\mathcal{T} = 4$, and $Ra=10^5$.}
    \label{fig:fig2_fields}
\end{figure}

Figure \ref{fig:fig2_fields} shows snapshots of the reference solution along with the difference with CDAnet-downscaled realizations obtained under different noise levels for the same time step.
Note that corresponding plots illustrating the temperature fields are presented in the Supplementary.
The spatial error distributions varies and increases in magnitude with increasing noise level, where the highest absolute errors are concentrated near the boundary layer and mixing fronts.
For small values of $\sigma_{obs}$, the difference between the downscaled and reference fields is negligible, with no discernible spatial pattern. 
As the noise level increases, the errors appear to concentrate near plumes and areas of intense mixing.
Moreover, the prediction presents an increasingly granular noise texture. 
For $\sigma_{obs}=0.1$, the large-scale structures are barely visible, and the prediction's reliability is lost for $\sigma_{obs} > 0.1$ with extremely noisy predictions.

\begin{figure}
    \centering
    \includegraphics[width=\linewidth]{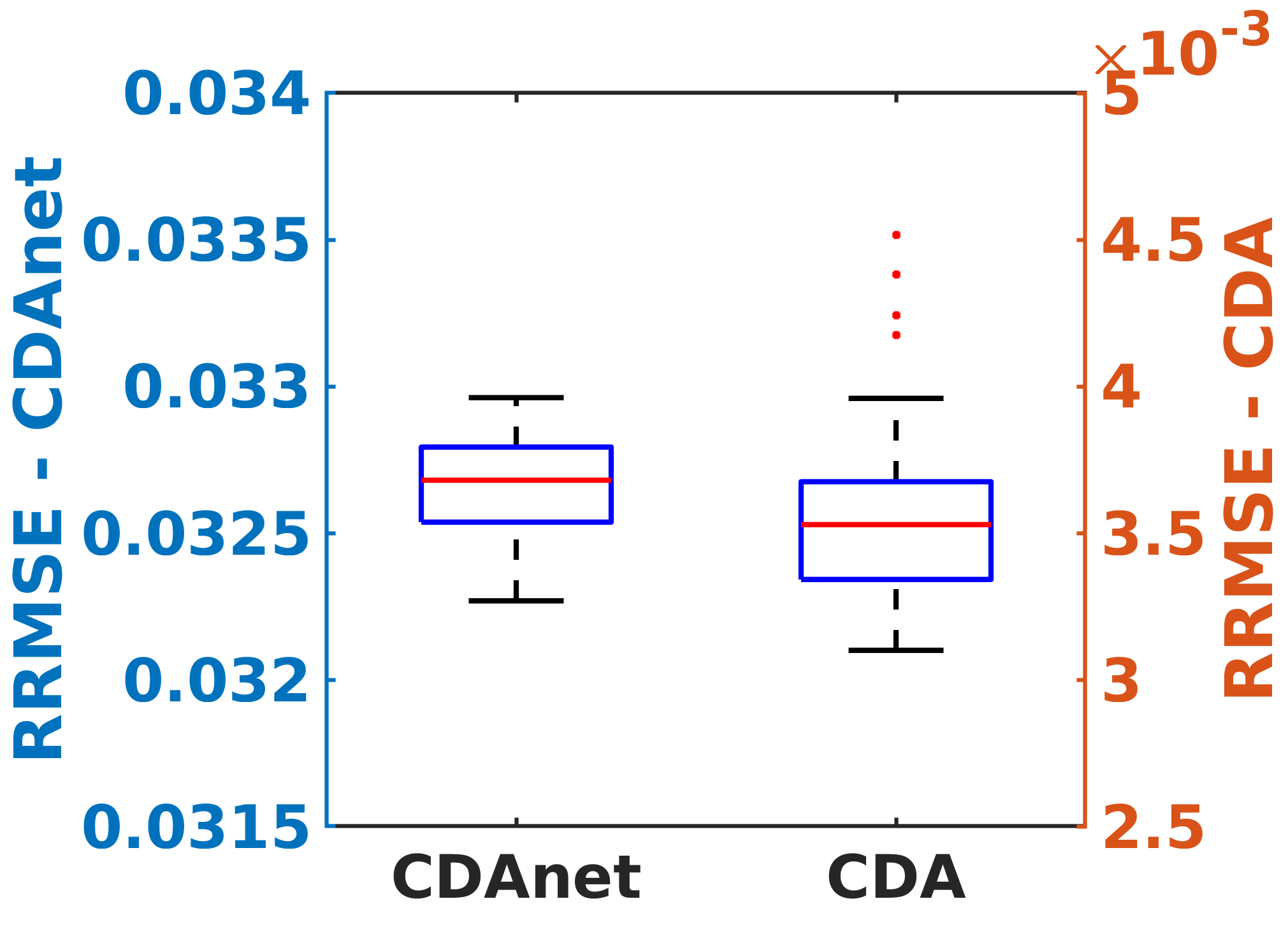}
    \caption{Boxplots of temperature RRMSE of CDAnet- and CDA-downscaled solutions. Each boxplot describes 50 realizations for $\sigma_{obs}=10^{-2}$. Both CDA and CDAnet were used to downscale 50 realizations of noisy coarse-scale solution trajectories with $\mathcal{S} = 4$, $\mathcal{T} = 4$, and $Ra=10^5$.}
    \label{fig:ensSize_CDAnetCDA}
\end{figure}

In the noise-free deterministic setting, CDA was shown to produce more reliable downscaled fields when compared to those obtained by CDAnet \citep{HammoudJAMES2022}. 
This result was attributed to CDAnet being a surrogate of the determining form map and hence has representation errors; i.e.\ limitations because it's an approximation of the determining form map. 
Consequently, CDA is considered as the benchmark when comparing downscaling algorithms.
Figure~\ref{fig:ensSize_CDAnetCDA} shows boxplots of the temperature RRMSE resulting from downscaling an ensemble of 50 coarse realizations using CDAnet and CDA with $\sigma_{obs}=0.01$.
The plot indicates that the RRMSEs of individual CDAnet-downscaled realizations are systematically larger than those of the CDA-downscaled ensemble. 
The CDAnet ensemble has an average RRMSE is approximately $3.28\%$, whereas that corresponding to the CDA ensemble is approximately $0.35\%$. 
This suggests that while the CDA downscaled fields are more accurate than those of CDAnet, it remains a dynamical downscaling algorithm requiring that the governing equations be solved on a high resolution mesh. 
Therefore, in the case when computational resources are limited, CDAnet remains a viable option for downscaling noisy coarse-scale realizations provided that an error overhead larger than that with CDA downscaled fields is associated with the CDAnet solution.

\begin{figure}[!htbp]
    \centering
    \subfloat[RMSE]{\includegraphics[width=60mm]{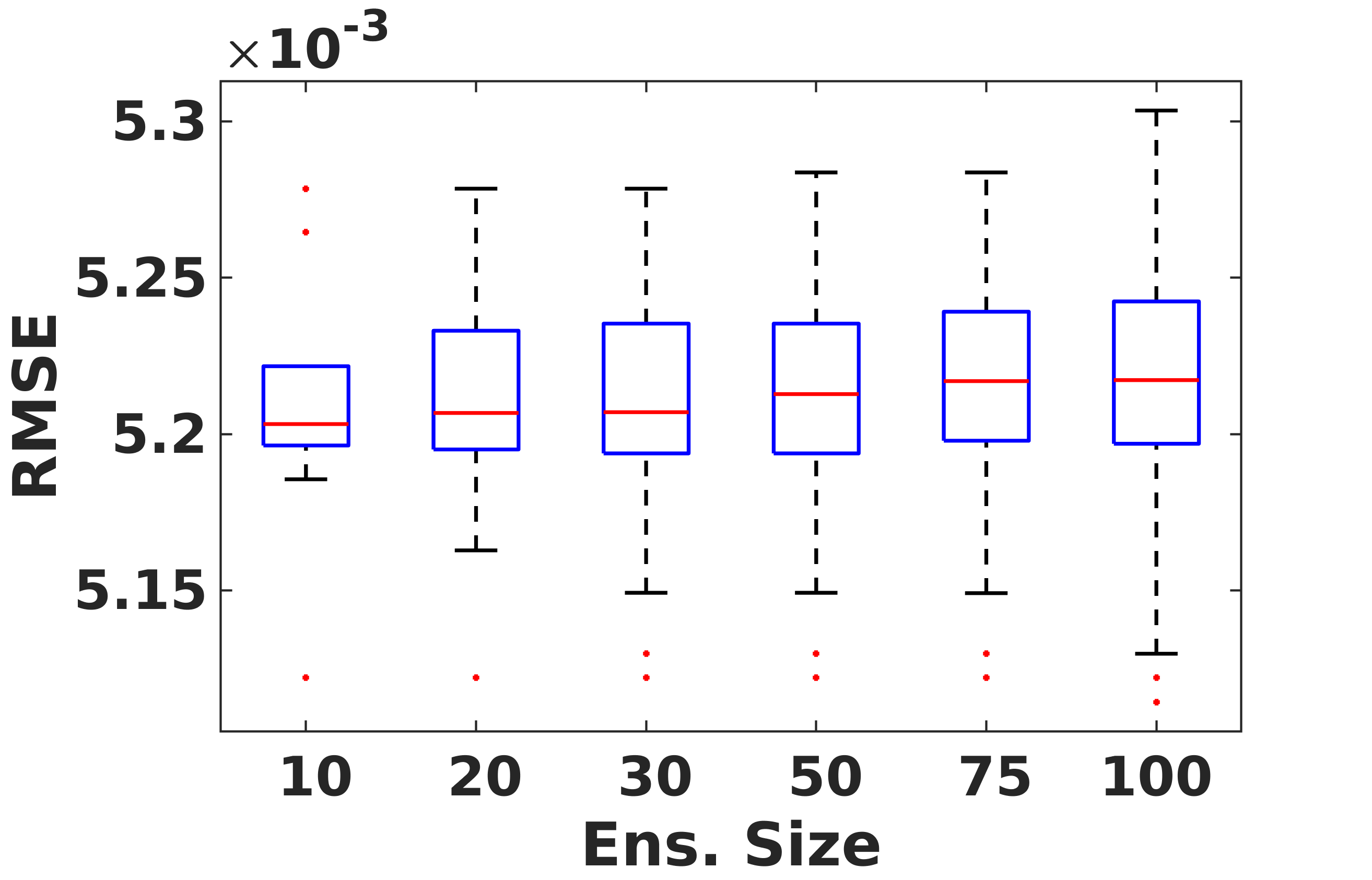}\label{part:aFIG1}} \,
    \subfloat[RRMSE]{\includegraphics[width=60mm]{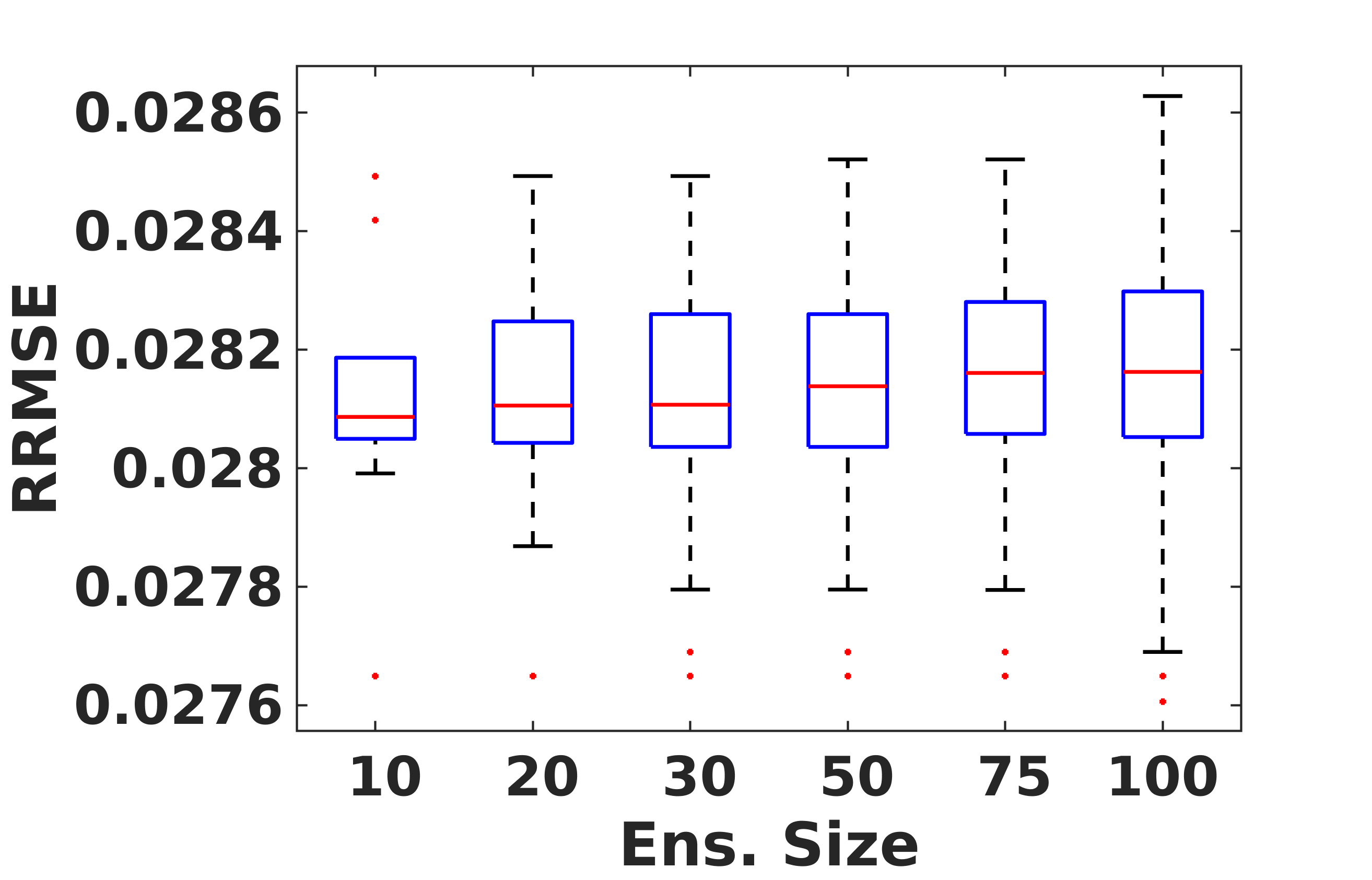}\label{part:bFIG1}} \,
\caption{Boxplots of the temperature (a) RMSE and (b) RRMSE of 50 CDAnet-downscaled fields for different ensemble sizes. In all cases, $\sigma_{obs}=0.01$, $\mathcal{S} = 4$, $\mathcal{T} = 4$, and $Ra=10^5$.}
\label{Fig:boxplots}
\end{figure}

The statistics of the downscaled ensembles are then examined by plotting in Figure \ref{Fig:boxplots} boxplots of the temperature RMSE and RRMSE for ensembles of increasing size.
The plots indicate that the distributions of the RMSE and RRMSE of the downscaled samples experience a negligibly small change for ensemble sizes of $20$ and larger, where the estimates of the mean and quantiles of the RMSE and RRMSE become almost identical for large ensembles.
This indicates that it is sufficient to consider moderate size ensembles to obtain suitable estimates of the errors' statistics associated with ensembles downscaled with CDAnet, at least as they are currently structured and trained.


The Kolmogorov-Smirnov test is applied to an ensemble of 50 CDAnet-downscaled realizations along the vertical centerlines and horizontal midline of $\Omega$, at different times for $\sigma_{obs}=0.01$, where the plots are presented in the supplementary for the sake of brevity.
The Kolmogorov-Smirnov is a statistical hypothesis test to examine whether samples are normally distributed.
For both time steps, the test indicates a p-value of one along all the spatial points along the centerline and midline, indicating that the downscaled fields are normally distributed along these directions. 
This suggests that CDAnet lifts a normally distributed ensemble of noisy coarse trajectories to a normally distributed ensemble of high-resolution trajectories preserving the distribution of the noisy realizations.

\begin{figure}
    \centering
    \includegraphics[width=\linewidth]{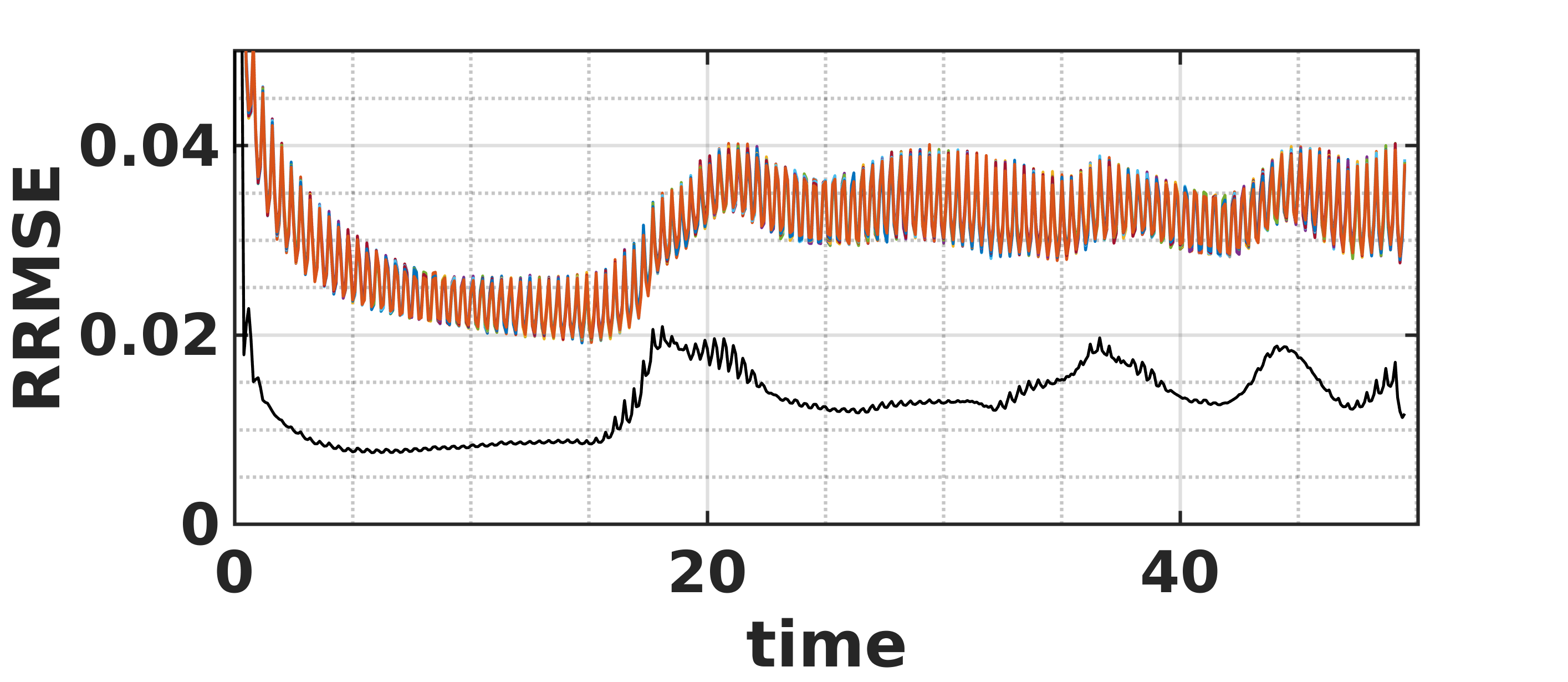}
    \caption{Curves showing the evolution of RRMSE of individual CDAnet-downscaled temperature distributions (colored) and the RRMSE of the corresponding ensemble average (black). The ensemble comprises 50 realizations, which were generated by downscaling an ensemble of coarse-scale solution trajectories with $\sigma_{obs}=0.01$, $\mathcal{S} = 4$, $\mathcal{T} = 4$, and $Ra=10^5$.}
    \label{fig:meanRRMSE_meanFieldRRMSE}
\end{figure}

Figure \ref{fig:meanRRMSE_meanFieldRRMSE} shows the individual temperature RRMSE curves corresponding to $50$ CDAnet-downscaled predictions along with the temperature RRMSE curves of the ensemble mean.  
The plot shows that RRMSE curves for the individual CDAnet-downscaled realizations decay over time before oscillating about an error plateau. 
The temperature RRMSE of the ensemble mean is appreciably lower than the RRMSE of the members of the downscaled ensemble, consistent with the analogue results of CDA \citep{HammoudCOMG2022}. 
These results suggest that by downscaling several noisy realizations using CDAnet and taking their average, an improved prediction can be obtained.

\subsubsection{Effect of Model Noise}
\label{res1_mod}

An ensemble of noisy CDAnet models is generated by sampling noise from a Gaussian distribution to perturb the parameters of each convolutional and dense layer of the CDAnet model during inference.
The noise-free low-resolution solution trajectory is input to the sampled noisy model realizations to predict their high-resolution counterparts. 
We consider model noise levels in a wide range, i.e., $\sigma_{mod} \in \left[0.001 ,\ 0.3 \right]$, with 50 realizations generated for each value of $\sigma_{mod}$. 

\begin{figure}
    \centering
    \includegraphics[width=\linewidth]{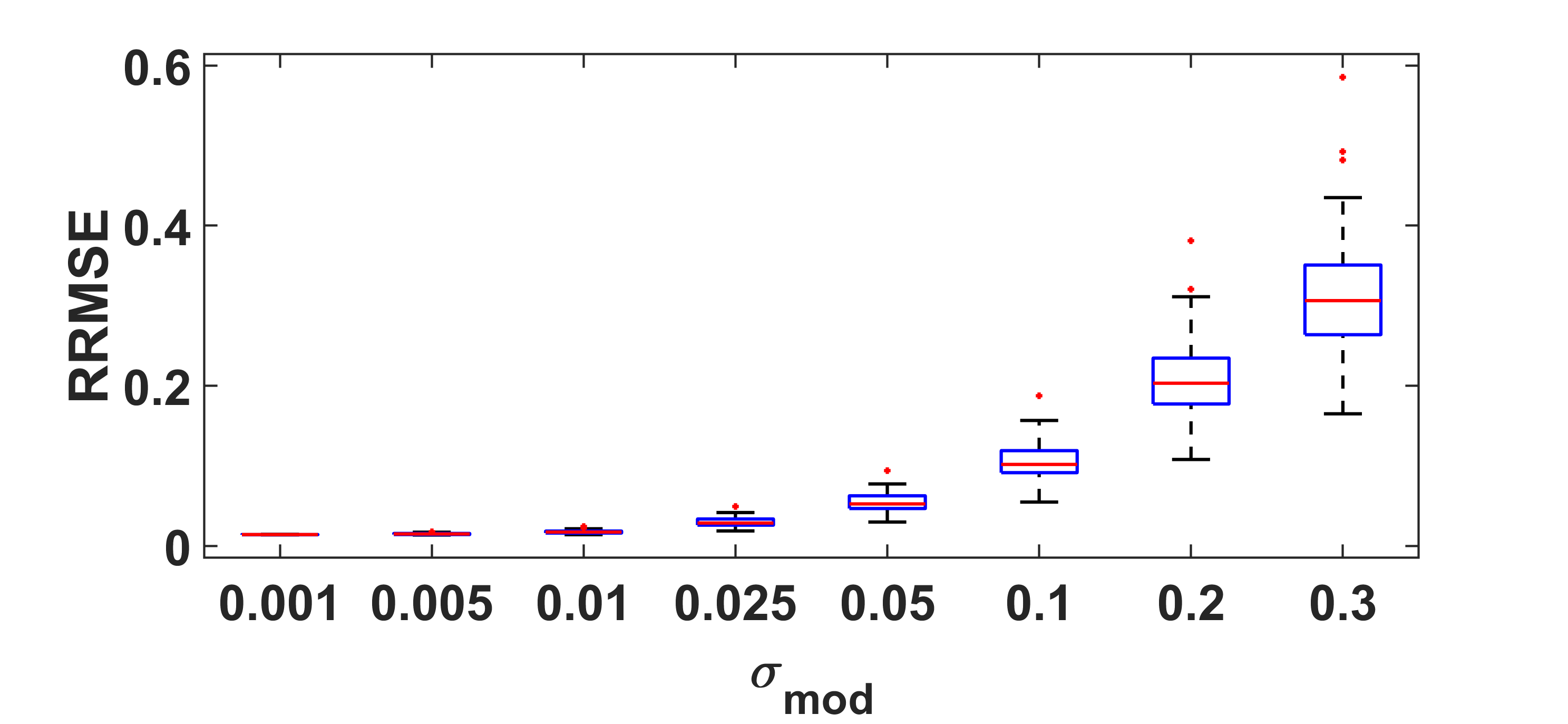}
    \caption{Boxplots of the temperature RRMSE of 50 CDAnet-downscaled fields for different values of $\sigma_{mod}$ with $\mathcal{S} = 4$, $\mathcal{T} = 4$, and $Ra=10^5$.}
    \label{fig:boxplots_RRMSE_modelErr}
\end{figure}

Figure \ref{fig:boxplots_RRMSE_modelErr} shows the boxplots of the distributions of the temperature RRMSE of 50 CDAnet-downscaled predictions. 
The RRMSE of the downscaled realizations increases with increasing $\sigma_{mod}$, where both the mean and standard deviation of the error's distribution increase with increasing $\sigma_{mod}$. 
This is, in part, different than the results obtained with observational noise only, where the errors resulting from different downscaled realizations collapse onto a single value, essentially having no spread at all.

\begin{figure}
    \centering
    \includegraphics[width=\linewidth]{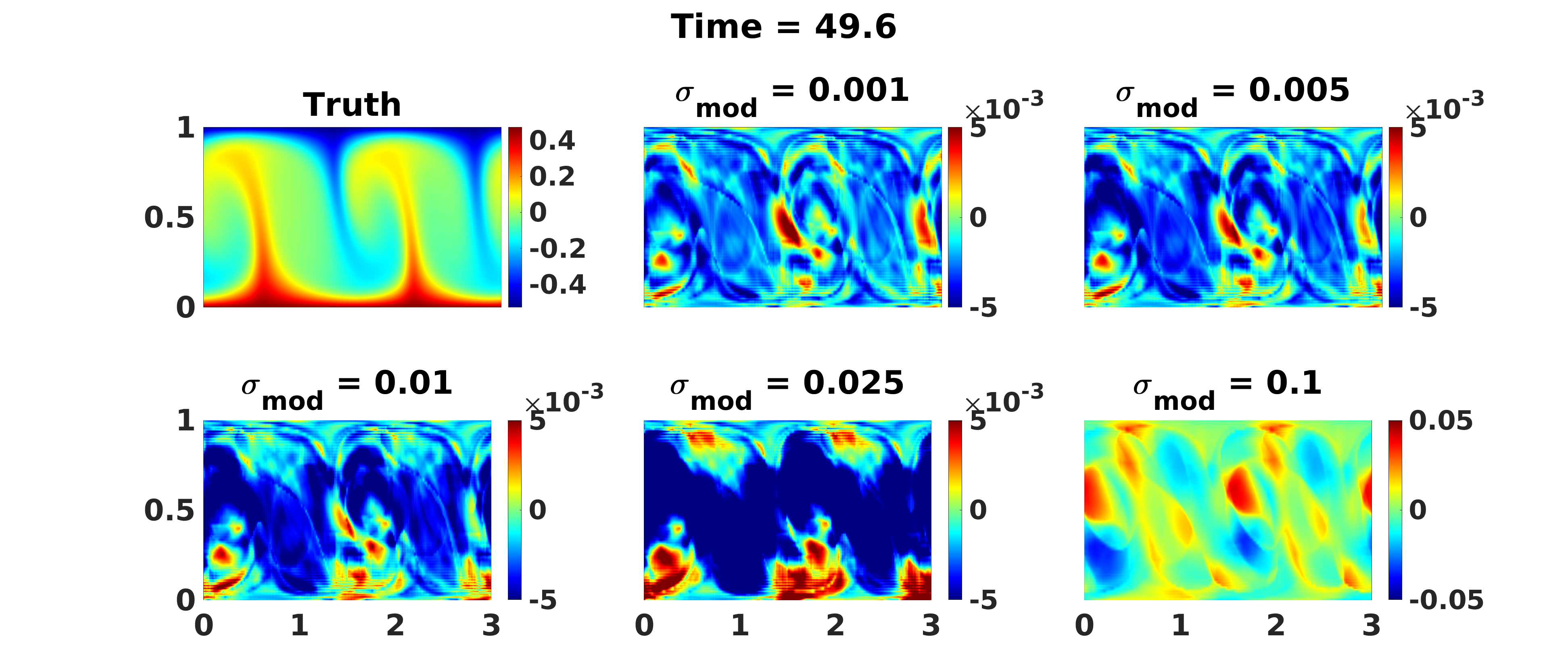}
    \caption{Snapshots of the reference temperature field along with the difference between a randomly selected sample downscaled temperature field and the high-resolution reference counterpart for different $\sigma_{mod}$ values at time $t=49.6$. Noisy CDAnet models were used to downscale a single noise-free coarse-scale solution trajectory with $\mathcal{S} = 4$, $\mathcal{T} = 4$, and $Ra=10^5$.}
    \label{fig:TFields_modelErr}
\end{figure}

Figure \ref{fig:TFields_modelErr} shows snapshots of the reference temperature solution along with its difference with CDAnet-downscaled realizations at $t=49.6$, for different $\sigma_{mod}$ as indicated.
Plots illustrating the temperature fields are presented in the Supplementary.
Unlike the case where only observational noise was considered, the error's spatial distribution does not correlate with the flow's dynamics.
In particular, the plot indicates that for $\sigma_{mod}\leq 0.01$, the differences with the reference solution are imperceivable, where the predicted large spatial scales of the downscaled solutions remain clear.
On the other hand, the finer-scales are distorted, where the prediction shows a slightly different position for the tip of the plume and with a different thickness.



\subsubsection{Effect of Combined Observational and Model Noises}
\label{res1_comb}

CDAnet models are then examined in the case when observational and model errors are present in the inference stage. 
Observational and model noise levels in the range $(\sigma_{obs},\sigma_{mod}) \in \left[0.001,\ 0.1\right] \times \left[0.001,\ 0.1\right]$ are considered.
For each $(\sigma_{obs},\sigma_{mod})$ pair selected, 50 CDAnet-downscaled solutions are generated for the analysis.


Table \ref{fig:combinedErrors_meanStd} outlines the means and standard deviations of the temperature RRMSE of 50 downscaled trajectories, generated for different ($\sigma_{obs}$, $\sigma_{mod}$) pairs. 
For $\sigma_{obs} < 0.01$ and $\sigma_{mod} < 0.025$, the mean RRMSE error falls below $4.5\%$, indicating that when data and model perturbations are small enough, CDAnet predictions are reliable. 
Furthermore, the average RRMSE increases appreciably with $\sigma_{obs}$ and at a faster rate in comparison to $\sigma_{mod}$; beyond $\sigma_{obs}=0.025$, CDAnet predictions become less reliable.
Unlike the mean errors, the standard deviation of the RRMSEs increase at faster rate with $\sigma_{mod}$, and exhibits a weaker dependence on $\sigma_{obs}$.
In other words, the variability in the errors associated with the downscaled fields are more sensitive to $\sigma_{mod}$ than to $\sigma_{obs}$, whereas their magnitude is more sensitive to $\sigma_{obs}$. 

To provide a qualitative representation of these results, snapshots of the downscaled temperature field for different $(\sigma_{obs}, \sigma_{mod})$ combinations are presented in the Supplementary Material.  
The plots show that by fixing $\sigma_{mod}$ and increasing $\sigma_{obs}$, the prediction becomes noisy with a granular texture covering the spatial domain.
In contrast when the $\sigma_{obs}$ is fixed, the imprint of the convolution filters appears stronger across the domain with increasing $\sigma_{mod}$.

The boxplots of the temperature RMSE and RRMSE for ensembles with different sizes, ranging from $10$ to $100$ with $\sigma_{obs} = 0.0125$ and $\sigma_{mod} = 0.01$ were plotted and analyzed.
These boxplots are shown in the supplementary because the conclusions drawn are similar to those discussed in the previous sections.
The results suggest that an ensemble size as small as 10 may be sufficient to capture the mean and spread of the error distribution, with both quantities exhibiting weak dependence on the ensemble size. 
This suggests that the CDAnet downscaled ensemble is rich in terms of the statistical information it provides, where the tails increase with increasing sample size, however, the mean estimate remains within 3$\%$ of that estimated by a small cardinality ensemble of downscaled fields. 
We also apply the Kolmogorov-Smirnov test to these downscaled realizations, in the same fashion as in Section \ref{res1_obs}.
In brief, the results showed that the downscaled samples are normally distributed, along the horizontal and vertical mid-sections of the domain.


The sensitivity of $\Lambda$ to different $Ra$ numbers is then examined in Figure \ref{fig:combinedErrors_effectOfRa}.
$\Lambda$ is computed using 50 downscaled realizations for each of $Ra=10^5$, $10^6$ and $10^7$.
Different CDAnet models that were comparably trained for different $Ra$ numbers were used for inference with $(\sigma_{obs}, \sigma_{mod}) = (0.01, 0.025)$ and $(\mathcal{S}, \mathcal{T})=(2, 2)$ to produce the realizations.
This configuration was selected to minimize the impact of CDAnet's representation errors and to assess its sensitivity to $Ra$ number.
Moreover, this configuration is the only common one between the $Ra$ numbers considered, as seen in \citep{HammoudJAMES2022}, where at $Ra=10^7$, the maximum $\mathcal{S}$ achievable is 2 for both CDA and CDAnet for the computational mesh considered.
The figure shows that $\Lambda$ slightly increases with increasing $Ra$ numbers, where the plots indicate a weak dependence on the $Ra$ number, provided it lies within the CDA guarantees for downscaling \citep{Azouani2013}.
Therefore, a trained CDAnet model is capable of downscaling flows at a wide range of $Ra$ numbers without compromising accuracy, which adheres to the theoretical guarantees of CDA and the determining form map.

Similarly, the sensitivity of $\Lambda$ to different $\mathcal{S}$ and $\mathcal{T}$ values was examined, as shown in Figure \ref{fig:combinedErrors_effectOfST}, by individually varying $\mathcal{S}$ and $\mathcal{T}$.
The CDAnet models were then used to generate 50 downscaled realizations with $(\sigma_{obs}, \sigma_{mod})=(0.005, 0.005)$ and $Ra=10^5$.
Note that the CDAnet model parameters were optimized for the specific $\mathcal{S}$ and $\mathcal{T}$ under consideration.
$\Lambda$ slightly increases with increased $\mathcal{S}$, suggesting that the average errors associated with an ensemble of CDAnet downscaled fields exhibits a weak dependence on $\mathcal{S}$, provided that CDAnet could accurately downscale noise-free coarse scale information at that resolution. 
This takeaway is different than the results obtained for CDA, where $\Lambda$ increases quadratically with $\mathcal{S}$ \citep{HammoudCOMG2022}.
On the other hand, a strong monotonic dependence of $\Lambda$ on $\mathcal{T}$ is observed, with $\Lambda$ increasing quadratically with $\mathcal{T}$, which represents a stronger relation in comparison to CDA's linear dependence on $\mathcal{T}$ \citep{HammoudCOMG2022}.



\subsection{Perturbed Training Conditions}
\label{sec:results2}

In a more general setting, training and validation data might be noisy requiring the model to be trained with uncertainties.
In this section, both coarse-scale data and model parameters are perturbed during training, such that the optimization process is affected by noise. 
Specifically, observational and model noise are sampled at the beginning of a training epoch and added to the low-resolution inputs and model parameters of CDAnet. 
The inputs were then propagated through CDAnet resulting in a higher-resolution prediction, which is used to optimize the network's parameters.
In particular, by perturbing the observational data and model weights, the learned model is hypothesized to become more robust, particularly less sensitive to perturbations \citep{Morcos2018_modelErrors, Tsai2021_modelErrors}.
The performance of CDAnet that is trained with observational and model perturbations for different combinations of $\left( \mathcal{S}, \mathcal{T}, Ra \right)$ is then assessed against different conditions of observational and model noise during inference.

Table \ref{table:noiseLevels_trainedSec2} lists the temperature RRMSE obtained during CDAnet's validation for different training observational and model noise levels.
Here, the considered noise levels correspond to the cases of no, low, medium, and high noise levels.
The results show that the average validation RRMSE of the CDAnet models trained using different noise levels are comparable, suggesting that the trained CDAnet model is weakly sensitive to observational and model perturbations within the considered range.
The results analyzed in Sections (\ref{res2_obs})-(\ref{res2_comb}) correspond to the CDAnet model that was trained with $(\sigma_{obs}, \sigma_{mod}) = (0.01, 0.01)$ (highlighted in bold font in Table \ref{table:noiseLevels_trainedSec2}).
This model was trained using the same procedure as the one presented by \citet{HammoudJAMES2022}, however, for the case presented here with $\mathcal{S} = \mathcal{T} = 4$ and $Ra=10^5$, the training parameters adopted consist of a learning rate of $0.1$, $\lambda = 0.01$ and a batch size of 20.
CDAnet is then employed in an inference setting to downscale a coarse-scale trajectory with different observational and model noise levels.

\begin{table}[!htbp]
\centering
\begin{tabular}{cc|cccc|}
\cline{3-6}
                                                      &       & \multicolumn{4}{c|}{$\sigma_{mod}$}                                                                                  \\ \cline{3-6} 
                                                      &       & \multicolumn{1}{c|}{0}           & \multicolumn{1}{c|}{0.001}       & \multicolumn{1}{c|}{0.01}        & 0.1         \\ \hline
\multicolumn{1}{|c|}{\multirow{4}{*}{$\sigma_{obs}$}} & 0     & \multicolumn{1}{c|}{0.01046} & \multicolumn{1}{c|}{0.0105} & \multicolumn{1}{c|}{0.0100} & 0.0112 \\ \cline{2-6} 
\multicolumn{1}{|c|}{}                                & 0.001 & \multicolumn{1}{c|}{0.00974} & \multicolumn{1}{c|}{0.0114}  & \multicolumn{1}{c|}{0.0124} & 0.0117 \\ \cline{2-6} 
\multicolumn{1}{|c|}{}                                & 0.01  & \multicolumn{1}{c|}{0.00985} & \multicolumn{1}{c|}{0.0115} & \multicolumn{1}{c|}{\textbf{0.0101}} & 0.0123 \\ \cline{2-6} 
\multicolumn{1}{|c|}{}                                & 0.1   & \multicolumn{1}{c|}{0.01198} & \multicolumn{1}{c|}{0.0097} & \multicolumn{1}{c|}{0.0099} & 0.0122 \\ \hline
\end{tabular}
\caption{Temperature RRMSE values obtained during the CDAnet validation for the cases when observational and model noises are present during training. The results are herein presented for $(\mathcal{S}, \mathcal{T}) = (4, 4)$ and $Ra=10^5$.}
\label{table:noiseLevels_trainedSec2}
\end{table}

\subsubsection{Effect of Observational Noise}
\label{res2_obs}

Similar to Section \ref{res1_obs}, the impact of $\sigma_{obs}$ on downscaling coarse-scale inputs was analyzed by adding noise to the inputs of the trained CDAnet without perturbing its parameters during inference. 
For each $\sigma_{obs}$, an ensemble of 50 downscaled realizations is generated using CDAnet using 50 coarse-scale noisy realizations of a reference solution trajectory that was not present in the training nor validation datasets.
Figure \ref{fig:sec2_obsNoise_1} presents boxplots of the temperature RRMSE for the downscaled solutions corresponding to different $\sigma_{obs}$.   
These results show that the estimates are fairly flat for $\sigma_{obs} \leq 0.01$; however, for larger $\sigma_{obs}$, the RRMSE of the downscaled fields rapidly increases with $\sigma_{obs}$.
Nevertheless, the variability remains minor across all $\sigma_{obs}$, where the standard deviation of the RRMSE of 50 downscaled realizations remains small for all $\sigma_{obs}$.
This suggests that the errors associated with CDAnet predictions exhibit weak variability when only observational errors are present, and similarly to the case of noise-free training conditions, the mean RRMSE increases, however all the downscaled realizations collapse about a similar RRMSE.


Figure \ref{fig:Sec1_2_lambda_sigObs} plots the curves for $\Lambda$ as a function of $\sigma_{obs}$ for the case of perfect and perturbed training conditions, in addition to  the curves of their theoretical dependence \citep{Bessaih2015}.
The results show that $\Lambda$ slowly increases for $\sigma_{obs}<0.01$ and $\sigma_{obs}<0.025$, in the cases of perfect and perturbed training conditions, respectively.
Furthermore, for $\sigma_{obs}<0.025$ and perfect training conditions, $\Lambda$ is smaller than that corresponding to the perturbed training conditions.
This behavior indicates that perturbations during training result in an additional source of uncertainty. 
Finally, $\Lambda$ increases quadratically with $\sigma_{obs}$ for $\sigma_{obs}>0.01$ and $\sigma_{obs}>0.025$ in the cases of perfect and perturbed training conditions, respectively.
This indicates that CDAnet and CDA experience similar dependency on $\sigma_{obs}$ \citep{Bessaih2015, HammoudCOMG2022}.

\subsubsection{Effect of Model Noise}
\label{res2_mod}

The impact of model noise on CDAnet downscaled fields is assessed by downscaling an unperturbed coarse-scale trajectory using different CDAnet model realizations, whose parameters are perturbed during training and inference.
This is equivalent to analyzing a perturbed version of the lifting function operating in a scenario where the coarse-scale information is obtained from an imperfect numerical model.
Figure \ref{fig:sec2_modNoise_1} presents the boxplots of the temperature RRMSE values of 50 CDAnet-downscaled fields for different values of $\sigma_{mod}$.
The figure shows that the mean RRMSE and its standard deviation increase with increasing $\sigma_{mod}$, where unlike the case involving only $\sigma_{obs}$, the ensemble spread here is considerable.

Figure \ref{fig:Sec1_2_lambda_sigMod} plots the curves for $\Lambda$ as a function of $\sigma_{mod}$ for the cases of perfect and perturbed training conditions as well as the best fit curve.
$\Lambda$ increases slowly with $\sigma_{mod}$ for $\sigma_{mod}<0.025$ and $\sigma_{mod}<0.075$ in the cases of perfect and perturbed training conditions, respectively.
For $\sigma_{mod}<0.075$, $\Lambda$ corresponding to perfect training conditions is smaller than that with perturbed training conditions.
This indicates that when model errors are present during training, additional uncertainties adversely impact the downscaling performance of CDAnet.
The plots indicate that $\Lambda$ increases quadratically with $\sigma_{mod}$ for $\sigma_{mod}>0.025$ and $\sigma_{mod}>0.075$ for the cases of perfect and perturbed training conditions, respectively.
This indicates that representation errors are prominent for a wide range of $\sigma_{mod}$.

\subsubsection{Effect of Combined Observational and Model Noises}
\label{res2_comb}

Table \ref{fig:sec2_aggNoise_meanStd} outlines the means and standard deviations of the RRMSE values of 50 CDAnet-downscaled temperature fields.
Both the mean and standard deviation of the temperature RRMSE increase with increasing $\sigma_{mod}$ and $\sigma_{obs}$.
Similar to Table \ref{fig:combinedErrors_meanStd}, the mean error is more sensitive to $\sigma_{obs}$ than to $\sigma_{mod}$.
On the other hand, the standard deviation of the error is more sensitive to $\sigma_{mod}$.
Furthermore, the values of the means and standard deviations of the RRMSE of 50 CDAnet-downscaled fields are larger than those obtained in the case of perfect training conditions (Table \ref{fig:combinedErrors_meanStd}).
Although training with observational and model noises yields CDAnet models with RRMSE errors $<2\%$, these models appear to be more sensitive to observational and model noises than those trained under perfect training conditions due to a compounding effect.


\subsection{Training with Imperfect Downscaled Data}
\label{sec:results3}

The performance of CDAnet models trained with imperfect datasets is investigated by first downscaling different coarse-scale solution trajectories using CDA with additional model errors incorporated to the downscaling. 
In particular, $30$ different solution trajectories of the RB convection with $Ra = 10^6$ are generated to serve as reference solutions train CDAnet models. 
CDA was then employed to downscale a coarse representation of each of these solution trajectories with $\mathcal{S}=4$ and $\mathcal{T} = 4$; however, when downscaling, $Ra$ was set to $1.3 \times 10^6$ resulting in imperfectly downscaled solutions.
During training, the PDE loss was computed with $Ra = 1.3\times 10^6$, i.e., the assumed Ra number.
The results are presented for $\mathcal{S} = \mathcal{T} = 4$, where CDAnet was trained with a learning rate of 0.15, $\lambda = 0.01$ and a batch size of 20 with the remaining training parameters held the same as in \citep{HammoudJAMES2022}.
The trained CDAnet models are then employed in the inference stage, in the presence of observational and model uncertainties, to downscale low-resolution solution trajectories of the RB convection at Ra numbers of $10^6$ and $1.3\times 10^6$, which are the reference and modified Ra numbers of the flow field. 
In this scenario, CDAnet was trained in the absence of observational and model noise.

\subsubsection{Effect of Observational Noise}
\label{res3_obs}

The CDAnet trained using the imperfectly downscaled dataset is first employed to downscale an ensemble of 50 noisy coarse-scale realizations of a reference RB solution at $Ra=10^6$ and $1.3\times 10^6$.
The RRMSE boxplots of temperature and velocity components are plotted for both Ra numbers and presented in the Supplementary Material. 
The boxplots indicate that the temperature and velocities RRMSEs both increase monotonically with $\sigma_{obs}$.
The RRMSE of the downscaled velocity fields are comparable with those obtained in the case when CDAnet was trained with uncertainties (Section \ref{sec:results2}).
On the other hand, the temperature RRMSE has a larger amplitude in comparison to those obtained in the previous two scenarios.
For instance, when $\sigma_{obs}\leq 0.01$, the temperature RRMSE values for both cases are approximately $10\%$, whereas the RRMSE of the velocity components is below $4\%$ for $Ra=10^6$ and approximately $13\%$ for $Ra=1.3\times10^6$.
The results suggest that when only observational noise are present during inference, training with imperfect downscaled fields leads to errors that are comparable to those obtained when training CDAnet in the presence of observational and model perturbations. 
This indicates that, for the case of observational noise, introducing observational and model noise during training results in a CDAnet model with comparable results to that trained with imperfectly downscaled data.

Figure \ref{fig:Sec3_lambda_sigObs} shows the variation of $\Lambda$ as a function of $\sigma_{obs}$ for $Ra$ numbers of $10^6$ and $1.3\times 10^6$.
$\Lambda$ increases monotonically with $\sigma_{obs}$ with a similar behavior for both $Ra$ numbers.
The plot shows that $\Lambda$ quadratically increases with $\sigma_{obs}$ for $\sigma_{obs}>0.05$, whereas for smaller $\sigma_{obs}$, the representation errors are larger, as indicated by the plateau in $\Lambda$.
These results indicate that CDAnet's training conditions impact the level of its representation errors but does not affect the behavior of $\Lambda$ as a function of $\sigma_{obs}$.
Hence, CDAnet is able to downscale coarse-scale observational data with a similar performance across different $Ra$ numbers in the deterministic noise-free setting \citep{HammoudJAMES2022} and in the uncertain setting in the presence of observational noise.

\subsubsection{Effect of Model Noise}
\label{res3_mod}

Figure \ref{fig:sec3_modNoise_1} illustrates boxplots of the RRMSE corresponding to an ensemble of 50 downscaled temperature and u-velocity fields on a log-scale. 
The RRMSE values for both variables increase monotonically as $\sigma_{mod}$ increases with a widening distribution as $\sigma_{mod}$ increases. 
The downscaled temperature fields exhibit larger RRMSE values compared to those of the velocity fields. 
Specifically, for $\sigma_{mod} \leq 0.05$, the RRMSEs of the temperature fields remain below $10\%$, whereas the RRMSE of the u-velocity is approximately half of those values for the same $\sigma_{mod}$. 
Importantly, the error values are comparable for both $Ra$ numbers, indicating the generalizability of the CDAnet model in downscaling coarse-solution trajectories belonging to nearby $Ra$ numbers \citep{HammoudJAMES2022}.
Here however, generalizability has two aspects, first the CDAnet's inference predictions result in comparable errors for solution trajectories with nearby $Ra$ numbers, and second, the resulting variance of the predictions are also comparable.


Figure \ref{fig:Sec3_lambda_sigMod} shows $\Lambda$ as a function of $\sigma_{mod}$ for the reference and assumed $Ra$ numbers. 
$\Lambda$ increases monotonically with $\sigma_{mod}$ and at a quadratic rate with $\sigma_{mod}$ for $\sigma_{mod}>0.1$.
In contrast, for smaller values of $\sigma_{mod}$, the curves reach a plateau, which is attributed to representation errors. 
These results indicate that the errors of a well-trained CDAnet model depend on $\sigma_{mod}^2$ if the effect of the added noise exceeds inherent representation errors.
These results further demonstrate that the CDAnet's performance generalizes for coarse-scale solution trajectories with different $Ra$ numbers even in the case when model errors are present.
Furthermore, CDAnet exhibits a similar dependence on the $\sigma_{mod}$ for all three training scenarios considered.


\subsubsection{Effect of Combined Observational and Model Noises}
\label{res3_comb}

Table \ref{fig:sec3_combinedNoise} and \ref{fig:sec3_combinedNoiseb} presents the means and standard deviations of the temperature and u-velocity RRMSE values of 50 CDAnet-downscaled fields for both $Ra$ numbers, respectively.
The mean RRMSE of both variables increases with increasing $\sigma_{obs}$ and $\sigma_{mod}$; however, the standard deviation increases at a faster rate with increasing $\sigma_{mod}$ than with $\sigma_{obs}$.
Although the mean temperature RRMSE values are comparable for the two downscaled trajectories, the errors corresponding to the solution trajectory with $Ra=10^6$ are slightly smaller than those corresponding to $Ra=1.3\times 10^6$.
The RRMSE values of the horizontal velocity are considerably smaller than those obtained for the temperature fields, indicating that the temperature field is more sensitive to perturbations than the velocities.
These results indicate that a CDAnet model trained with downscaled fields, and trained with a $Ra$ number not far from the true one, is able to generalize its performance to nearby $Ra$ numbers achieving acceptable downscaling errors when both observational and model errors are involved.
Furthermore, these results indicate that CDAnet behaves similarly against $\sigma_{obs}$ and $\sigma_{mod}$ regardless of the training conditions, however, the representation errors become larger the farther the training conditions are from ideal underscoring the importance of the quality of data when employing PI-DNNs for downscaling.



\section{Conclusion}
\label{sec:conclusion}

In this work, we investigate the performance of CDAnet within an uncertain framework involving observational and model noise by examining three different training scenarios pertaining to idealized conditions, perturbed conditions, and training with imperfectly downscaled data. 
CDAnet is a physics informed deep neural network that acts as a surrogate of the determining form map, which is a theoretical mapping from the coarse data to its higher resolution counterpart.
The Rayleigh-B\'enard convection paradigm was adopted to test the performance of CDAnet in these scenarios.

In the first setting, CDAnet models were trained without noise resulting in the ``best'' surrogate of the determining form map, which enables a thorough comparison with theoretical guarantees from the literature \citep{Bessaih2015}.
During inference, observational and model noise were incorporated to generate different downscaled solution samples.
The errors associated with the downscaled fields of all solution variables increase with increasing observational ($\sigma_{obs}$) and model ($\sigma_{mod}$) noise levels.
However, the standard deviation of these errors is negligible when only $\sigma_{obs}$ is considered.
In contrast, when only $\sigma_{mod}$ is considered, the standard deviation of the errors increases with increasing $\sigma_{mod}$.
When both observational and model noises are considered, the mean errors of the downscaled ensemble is more sensitive to $\sigma_{obs}$ than to $\sigma_{mod}$.
Furthermore, the ensemble standard deviation is more sensitive to the increase in the $\sigma_{mod}$ than in $\sigma_{obs}$ level distribution.

Within this idealized training setting, CDAnet was shown to yield normally-distributed downscaled fields when the input noise is sampled from a Gaussian distribution, behaving similarly to CDA.
The statistical analysis also shows that the distribution of the downscaled fields' errors are represented efficiently with small ensembles, where beyond 20 realizations, changes in the errors' distribution are negligible with increasing sample size. 
The results also suggest that $\Lambda$ increases quadratically with $\mathcal{T}$, but has a much weaker sensitivity to $\mathcal{S}$ and $Ra$, provided selection of these parameters satisfies conditions for convergence of downscaling or existence of the lifting function \citep{Foias2012, Azouani2013, Foias_2014}. 
The dependence of $\Lambda$ on $\mathcal{S}$ and $\mathcal{T}$ is different than that observed with CDA \citep{HammoudCOMG2022}, which indicates a quadratic and linear dependence, respectively.
Finally, the RRMSE of the ensemble averaged downscaled field is always smaller than the RRMSE of the individual downscaled solution fields, suggesting that it is worthwhile employing CDAnet to generate a better estimate of the downscaled solution by sampling noisy trajectories, as is the case with CDA \citep{HammoudCOMG2022}.


In the second scenario, observational and model noise were incorporated during CDAnet's training and inference.
The RRMSE of the downscaled fields increases with increasing noise level, and that CDAnet displays a similar performance to that exhibited for perfect training conditions with an additional error. 
In particular, when only $\sigma_{obs}$ is considered, the RRMSE of all the different downscaled realizations are identical to each other with a small spread similarly to CDA, whereas the spread increases appreciably when $\sigma_{mod}$ is considered.
$\Lambda$ exhibits a plateau for small $\sigma_{obs}$ and $\sigma_{mod}$ due to representation errors, but behaves quadratically with noise level beyond a certain threshold.
Furthermore, when both $\sigma_{obs}$ and $\sigma_{mod}$ are considered, the spread of the downscaled RRMSE is more sensitive to the $\sigma_{mod}$ than to the $\sigma_{obs}$, and the reversed relationship holds for the mean RRMSE. 
Finally, the RRMSE values are generally larger in this setting than those obtained under perfect training conditions, suggesting that the network is insensitive to small perturbations that are smaller than representation errors. 


Finally, CDAnet models were trained using imperfect data obtained by downscaling coarse-scale solution trajectories using CDA with an incorrect $Ra$ number specified when downscaling; i.e.\ setting the incorrect viscosity to the downscaling system.
The trained CDAnet model was then employed to downscale coarse-scale solution trajectories of the Rayleigh-B\'enard corresponding to the true and perturbed $Ra$ numbers. 
The errors are comparable for all the considered combinations of $\sigma_{obs}$ and $\sigma_{mod}$, suggesting the generalizability of a well-trained CDAnet model at downscaling flows with different $Ra$ numbers in an uncertain framework.
Moreover, the errors exhibit a similar behavior in terms of the rate of increase, spread, and functional dependence on the noise level as the two previous cases. 
However, the error values are larger because of the combined effect of the imperfectly downscaled training data, imperfect PDE loss, and perturbations involved during inference. 

The results of this study demonstrate the effectiveness of CDAnet at downscaling coarse-scale trajectories in an uncertain setting, inclusive of model and observational noise. 
Future work includes extending the current framework to downscaling statistical distributions as opposed to individual realizations.
In terms of practical applications, we will focus on downscale realistic large-scale ocean and atmospheric fields.



\section*{Acknowledgements}
Research reported in this publication was supported by the Office of Sponsored Research (OSR) at King Abdullah University of Science and Technology (KAUST) CRG Award \#CRG2020-4336 and Virtual Red Sea Initiative Award \#REP/1/3268-01-01. The work of E.S.T. was supported in part by NPRP grant \# S-0207-200290 from the Qatar National Research Fund (a member of Qatar Foundation).


\clearpage
\newpage

\bibliographystyle{elsarticle-harv}
\bibliography{references.bib}

\clearpage
\newpage



\clearpage
\newpage


\clearpage
\newpage

\clearpage
\newpage

\clearpage
\newpage

\begin{figure}
    \centering
    \includegraphics[width=\linewidth]{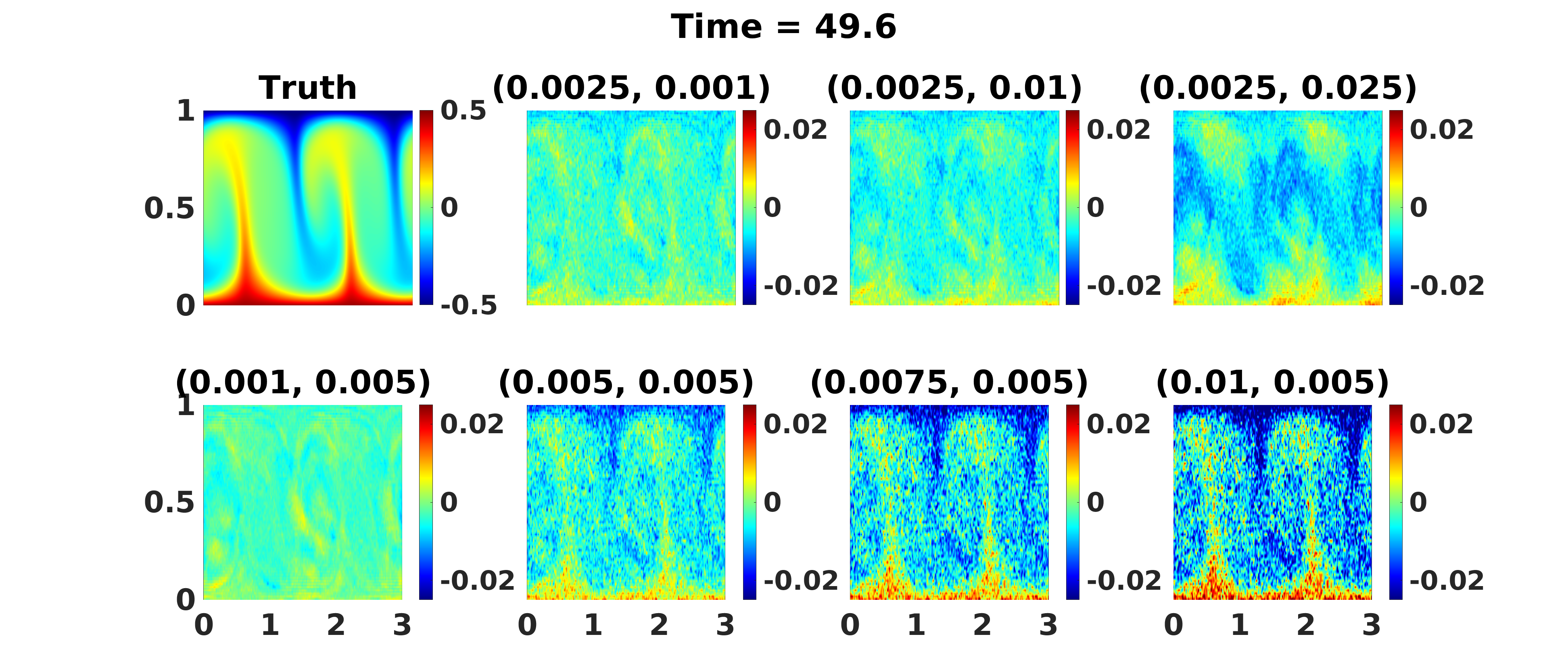}
    \caption{Snapshots of the reference temperature field along with the difference between the CDAnet-downscaled temperature field and reference counterpart for different $(\sigma_{obs}, \sigma_{mod})$ pairs. Noisy CDAnet models were used to downscale an ensemble of noisy representation of a coarse-scale reference trajectory with $\mathcal{S} = 4$, $\mathcal{T} = 4$, and $Ra=10^5$.}
    \label{fig:fields_combinedErrors}
\end{figure}


\begin{table}[!htbp]
\centering
\begin{center}
\resizebox{\textwidth}{!}{%
\begin{tabular}{|cc|cccccc|}
\hline
\multicolumn{2}{|c|}{\multirow{2}{*}{}} &
  \multicolumn{6}{c|}{$\sigma_{mod}$} \\ \cline{3-8} 
\multicolumn{2}{|c|}{} &
  \multicolumn{1}{c|}{0.001} &
  \multicolumn{1}{c|}{0.005} &
  \multicolumn{1}{c|}{0.01} &
  \multicolumn{1}{c|}{0.025} &
  \multicolumn{1}{c|}{0.05} &
  0.1 \\ \hline
\multicolumn{1}{|c|}{\multirow{8}{*}{$\sigma_{obs}$}} &
  0.001 &
  \multicolumn{1}{c|}{0.014303290} &
  \multicolumn{1}{c|}{0.015206472} &
  \multicolumn{1}{c|}{0.017790065} &
  \multicolumn{1}{c|}{0.030158020} &
  \multicolumn{1}{c|}{0.054953226} &
  0.106880534 \\ \cline{2-8} 
\multicolumn{1}{|c|}{} &
  0.0025 &
  \multicolumn{1}{c|}{0.015844652} &
  \multicolumn{1}{c|}{0.016664801} &
  \multicolumn{1}{c|}{0.019054420} &
  \multicolumn{1}{c|}{0.030930742} &
  \multicolumn{1}{c|}{0.055383338} &
  0.107099146 \\ \cline{2-8} 
\multicolumn{1}{|c|}{} &
  0.005 &
  \multicolumn{1}{c|}{0.020339128} &
  \multicolumn{1}{c|}{0.020987622} &
  \multicolumn{1}{c|}{0.022944854} &
  \multicolumn{1}{c|}{0.033514069} &
  \multicolumn{1}{c|}{0.056883720} &
  0.107874586 \\ \cline{2-8} 
\multicolumn{1}{|c|}{} &
  0.0075 &
  \multicolumn{1}{c|}{0.026100557} &
  \multicolumn{1}{c|}{0.026611799} &
  \multicolumn{1}{c|}{0.028194535} &
  \multicolumn{1}{c|}{0.037369142} &
  \multicolumn{1}{c|}{0.059278722} &
  0.109151615 \\ \cline{2-8} 
\multicolumn{1}{|c|}{} &
  0.01 &
  \multicolumn{1}{c|}{0.032465007} &
  \multicolumn{1}{c|}{0.032878419} &
  \multicolumn{1}{c|}{0.034180114} &
  \multicolumn{1}{c|}{0.042126072} &
  \multicolumn{1}{c|}{0.062447359} &
  0.11090924 \\ \cline{2-8} 
\multicolumn{1}{|c|}{} &
  0.025 &
  \multicolumn{1}{c|}{0.074234259} &
  \multicolumn{1}{c|}{0.074409680} &
  \multicolumn{1}{c|}{0.074991312} &
  \multicolumn{1}{c|}{0.078993776} &
  \multicolumn{1}{c|}{0.091704671} &
  0.129870645 \\ \cline{2-8} 
\multicolumn{1}{|c|}{} &
  0.05 &
  \multicolumn{1}{c|}{0.146174941} &
  \multicolumn{1}{c|}{0.146248200} &
  \multicolumn{1}{c|}{0.146524344} &
  \multicolumn{1}{c|}{0.148554830} &
  \multicolumn{1}{c|}{0.155675900} &
  0.181036527 \\ \cline{2-8} 
\multicolumn{1}{|c|}{} &
  0.1 &
  \multicolumn{1}{c|}{0.290167078} &
  \multicolumn{1}{c|}{0.290170259} &
  \multicolumn{1}{c|}{0.290268487} &
  \multicolumn{1}{c|}{0.291178685} &
  \multicolumn{1}{c|}{0.294711592} &
  0.308808656 \\ \hline
\end{tabular}%
}
\caption*{(a) Mean RRMSE}
\end{center}

\begin{center}
\resizebox{\textwidth}{!}{%
\begin{tabular}{|cc|llllll|}
\hline
\multicolumn{2}{|c|}{\multirow{2}{*}{}} &
  \multicolumn{6}{c|}{$\sigma_{mod}$} \\ \cline{3-8} 
\multicolumn{2}{|c|}{} &
  \multicolumn{1}{c|}{0.001} &
  \multicolumn{1}{c|}{0.005} &
  \multicolumn{1}{c|}{0.01} &
  \multicolumn{1}{c|}{0.025} &
  \multicolumn{1}{c|}{0.05} &
  \multicolumn{1}{c|}{0.1} \\ \hline
\multicolumn{1}{|c|}{\multirow{8}{*}{$\sigma_{obs}$}} &
  0.001 &
  \multicolumn{1}{l|}{0.000194} &
  \multicolumn{1}{l|}{0.001003} &
  \multicolumn{1}{l|}{0.002148} &
  \multicolumn{1}{l|}{0.006019} &
  \multicolumn{1}{l|}{0.012532} &
  0.025527 \\ \cline{2-8} 
\multicolumn{1}{|c|}{} &
  0.0025 &
  \multicolumn{1}{l|}{0.000176} &
  \multicolumn{1}{l|}{0.000917} &
  \multicolumn{1}{l|}{0.002015} &
  \multicolumn{1}{l|}{0.005881} &
  \multicolumn{1}{l|}{0.012441} &
  0.025478 \\ \cline{2-8} 
\multicolumn{1}{|c|}{} &
  0.005 &
  \multicolumn{1}{l|}{0.000141} &
  \multicolumn{1}{l|}{0.000736} &
  \multicolumn{1}{l|}{0.001699} &
  \multicolumn{1}{l|}{0.005468} &
  \multicolumn{1}{l|}{0.012136} &
  0.025305 \\ \cline{2-8} 
\multicolumn{1}{|c|}{} &
  0.0075 &
  \multicolumn{1}{l|}{0.000115} &
  \multicolumn{1}{l|}{0.000592} &
  \multicolumn{1}{l|}{0.001410} &
  \multicolumn{1}{l|}{0.004961} &
  \multicolumn{1}{l|}{0.011686} &
  0.025027 \\ \cline{2-8} 
\multicolumn{1}{|c|}{} &
  0.01 &
  \multicolumn{1}{l|}{0.0000996} &
  \multicolumn{1}{l|}{0.000491} &
  \multicolumn{1}{l|}{0.001186} &
  \multicolumn{1}{l|}{0.004462} &
  \multicolumn{1}{l|}{0.011150} &
  0.024658 \\ \cline{2-8} 
\multicolumn{1}{|c|}{} &
  0.025 &
  \multicolumn{1}{l|}{0.000085} &
  \multicolumn{1}{l|}{0.000281} &
  \multicolumn{1}{l|}{0.000637} &
  \multicolumn{1}{l|}{0.002585} &
  \multicolumn{1}{l|}{0.007968} &
  0.021415 \\ \cline{2-8} 
\multicolumn{1}{|c|}{} &
  0.05 &
  \multicolumn{1}{l|}{0.000131} &
  \multicolumn{1}{l|}{0.000302} &
  \multicolumn{1}{l|}{0.000583} &
  \multicolumn{1}{l|}{0.001802} &
  \multicolumn{1}{l|}{0.005305} &
  0.016314 \\ \cline{2-8} 
\multicolumn{1}{|c|}{} &
  0.1 &
  \multicolumn{1}{l|}{0.000241} &
  \multicolumn{1}{l|}{0.000503} &
  \multicolumn{1}{l|}{0.000920} &
  \multicolumn{1}{l|}{0.002287} &
  \multicolumn{1}{l|}{0.004991} &
  0.012702 \\ \hline
\end{tabular}%
}
\caption*{(b) Standard Deviation}
\end{center}

 \caption{Temperature RRMSE: (a) mean and (b) standard deviation for different combinations of $\sigma_{obs}$ and $\sigma_{mod}$. Noisy CDAnet models were used to downscale an ensemble of noisy representation of a coarse-scale reference trajectory with $\mathcal{S} = 4$, $\mathcal{T} = 4$, and $Ra=10^5$.
 \label{fig:combinedErrors_meanStd}}
\end{table}


\clearpage
\newpage

\begin{figure}
    \centering
    \includegraphics[width=\linewidth]{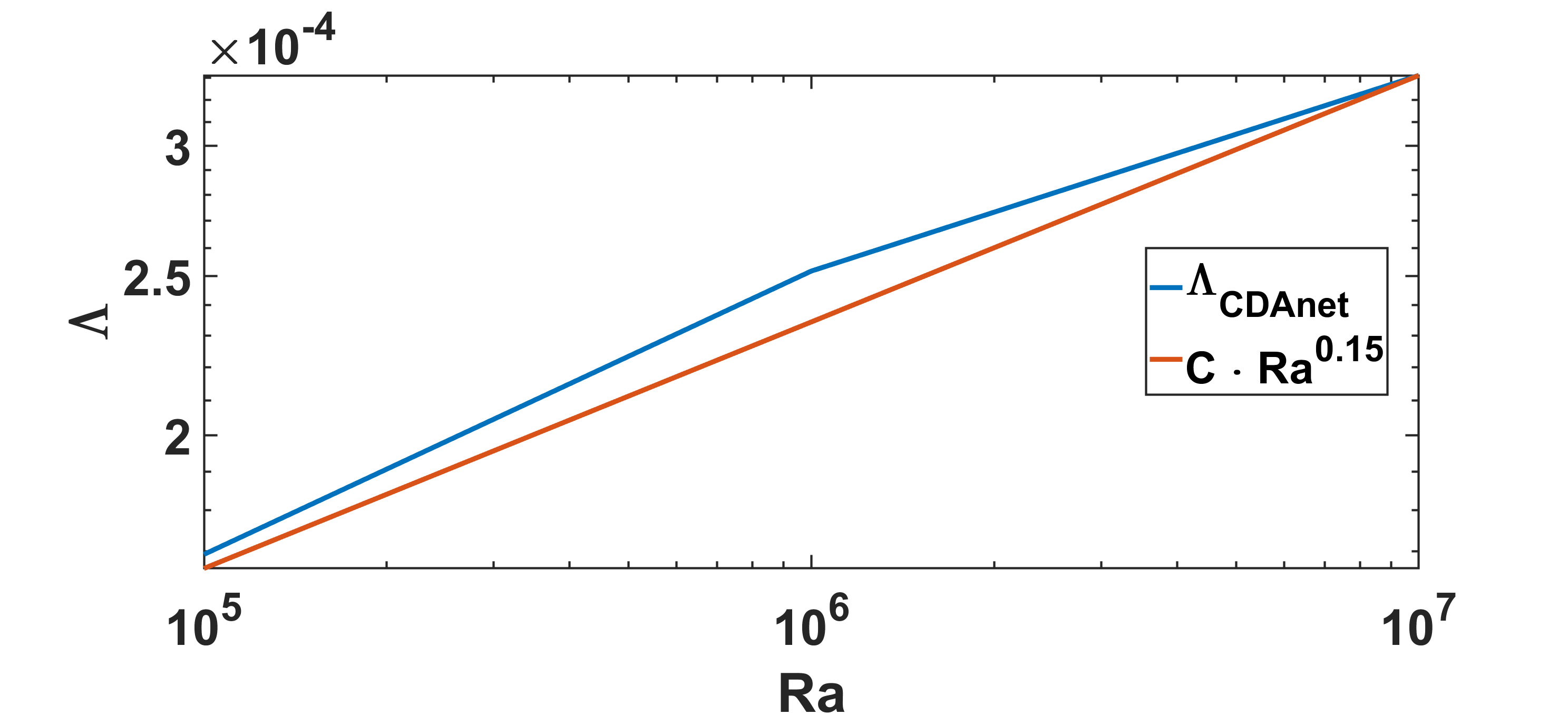}
    \caption{Plot representing the variations in $\Lambda$ with $Ra$. The values of $\Lambda$ are computed using 50 CDAnet predictions of the temperature fields. $(\mathcal{S}, \mathcal{T}) = (2, 2)$ and $(\sigma_{obs}, \sigma_{mod}) = (0.01, 0.025)$. The choice of $\mathcal{S}$ and $\mathcal{T}$ was made to ensure that representation errors are minimal.}
    \label{fig:combinedErrors_effectOfRa}
\end{figure}

\begin{figure}
    \centering
    \subfloat[]{\includegraphics[width=\linewidth]{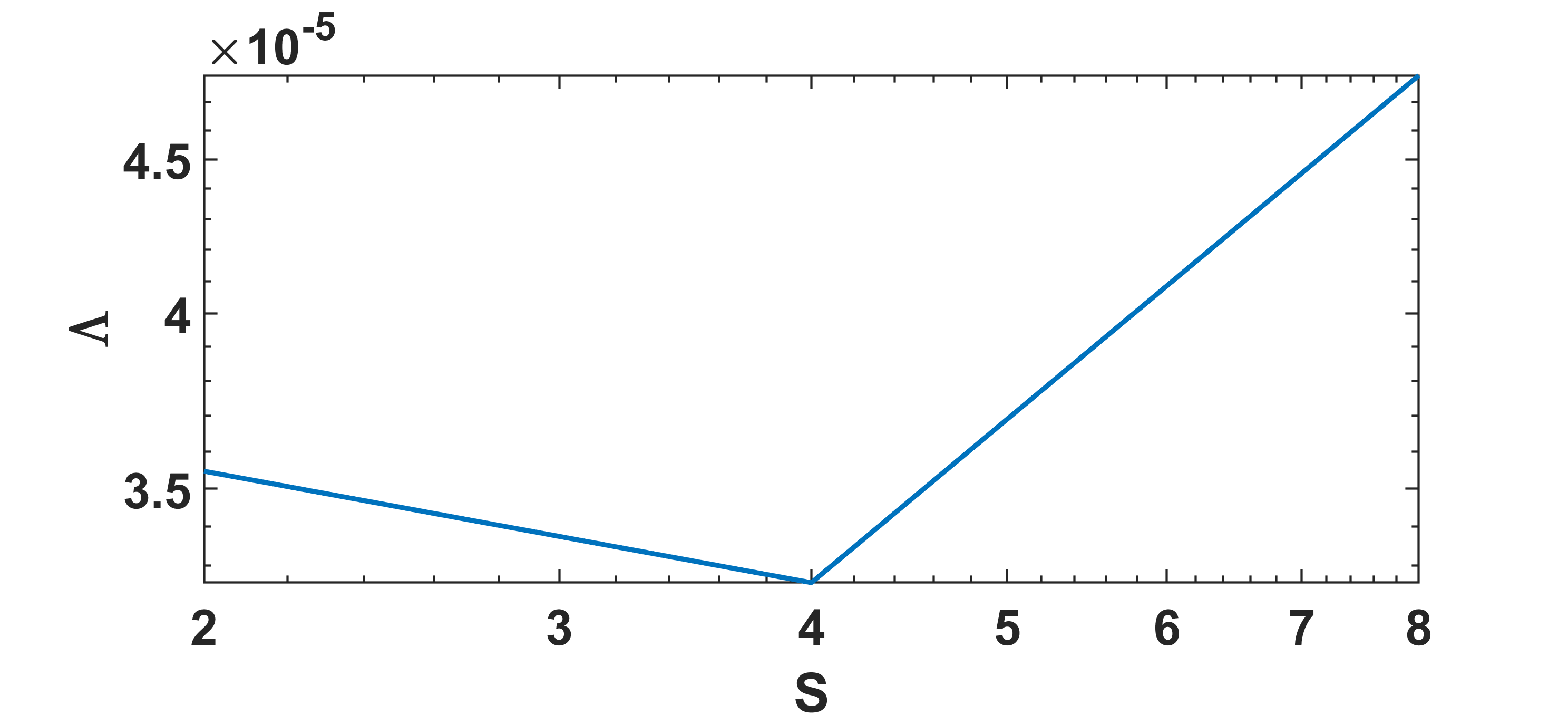}} \,
    \subfloat[]{\includegraphics[width=\linewidth]{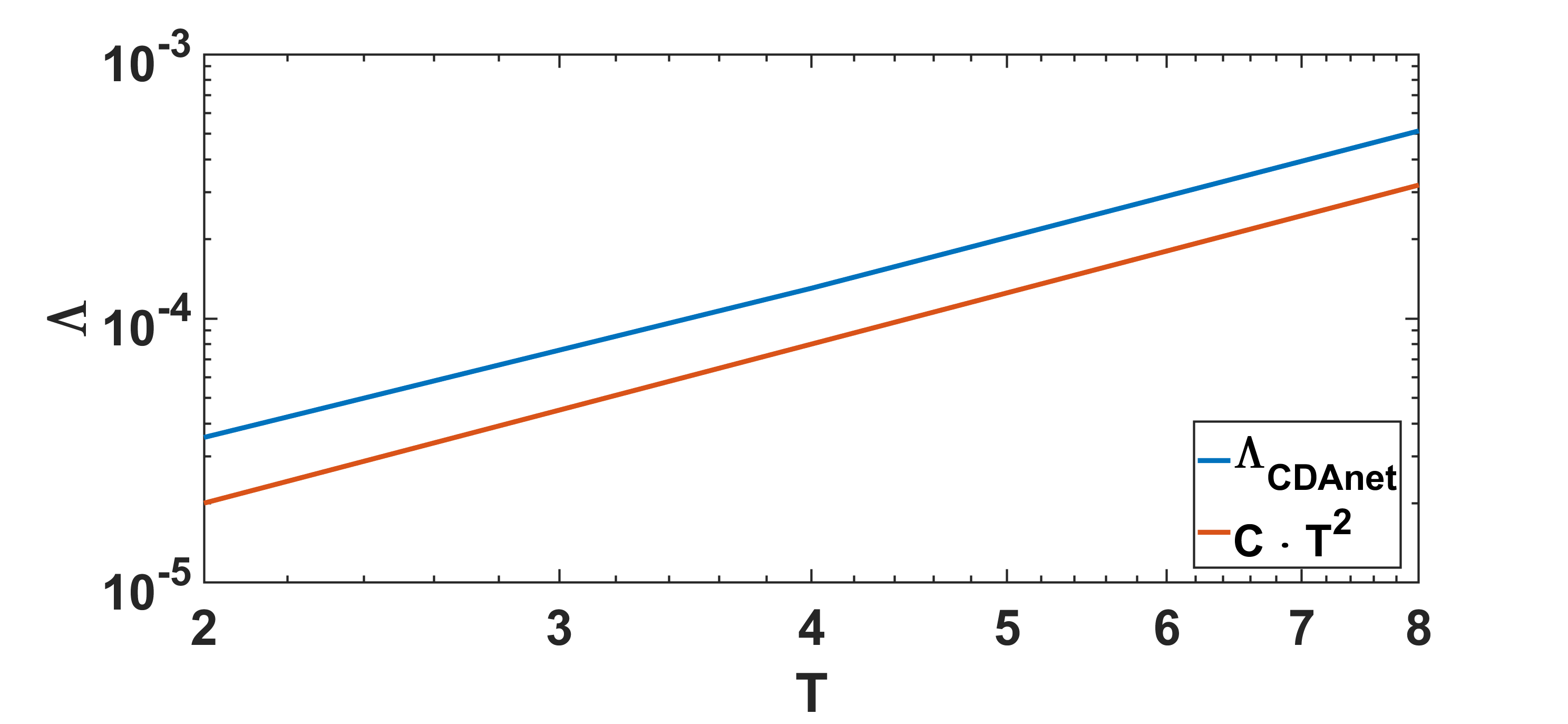}} \,
    \caption{Curves of $\Lambda$ as a function of (a) $\mathcal{S}$ and (b) $\mathcal{T}$ for $Ra = 10^5$, for $(\sigma_{obs}, \sigma_{mod})=(0.005, 0.005)$. In all cases, $\Lambda$ is computed for 50 realizations of the downscaled solution.}
    \label{fig:combinedErrors_effectOfST}
\end{figure}




\clearpage
\newpage


\begin{figure}[!htbp]
    \centering
    \includegraphics[width=\linewidth]{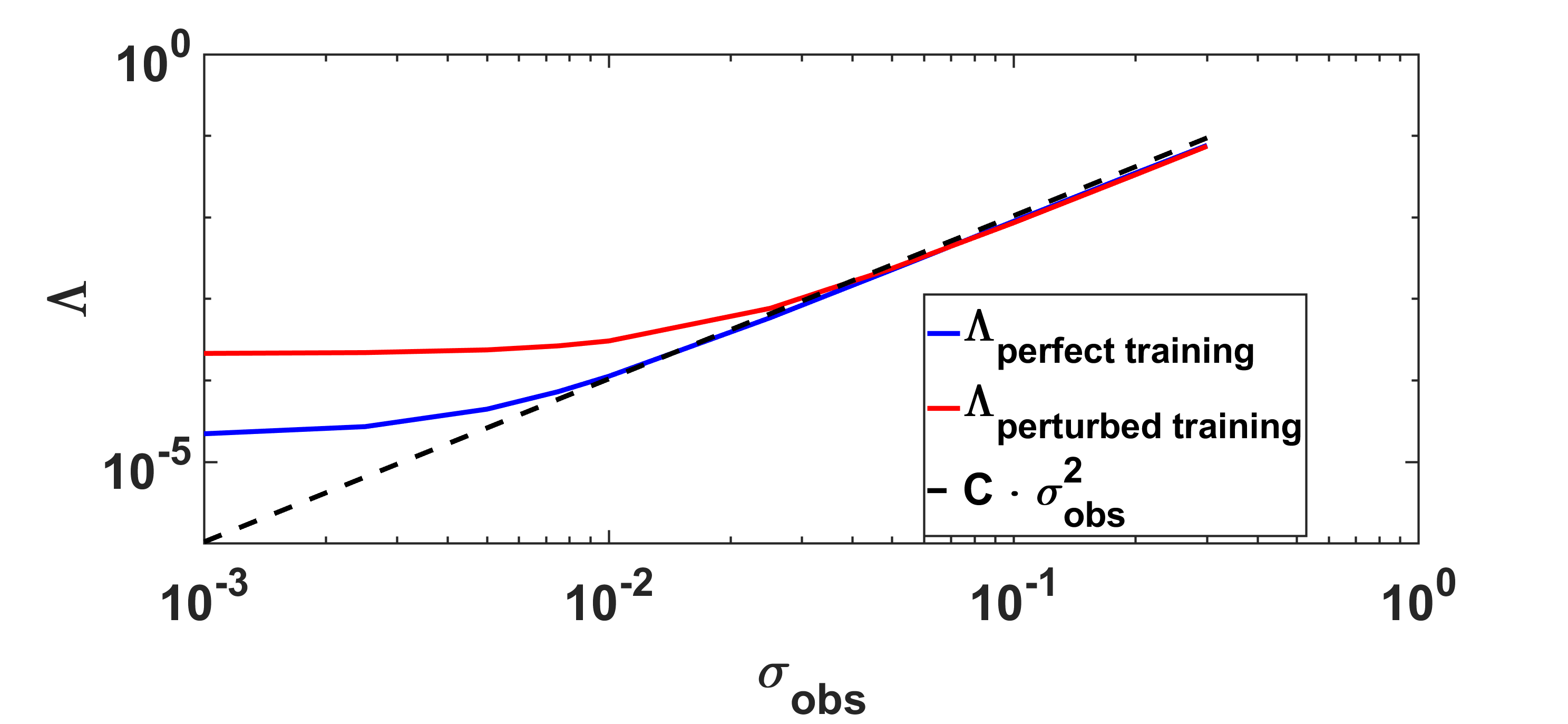}
    \caption{Curves of $\Lambda$ as a function of $\sigma_{obs}$ are shown for different $\sigma_{obs}$ when CDAnet was trained under perfect and perturbed conditions; $(\mathcal{S}, \mathcal{T}) = (4, 4)$ and $Ra=10^5$. Note that C represents a positive constant.}
    \label{fig:Sec1_2_lambda_sigObs}
\end{figure}

\begin{figure}[!htbp]
    \centering
    \includegraphics[width=150mm]{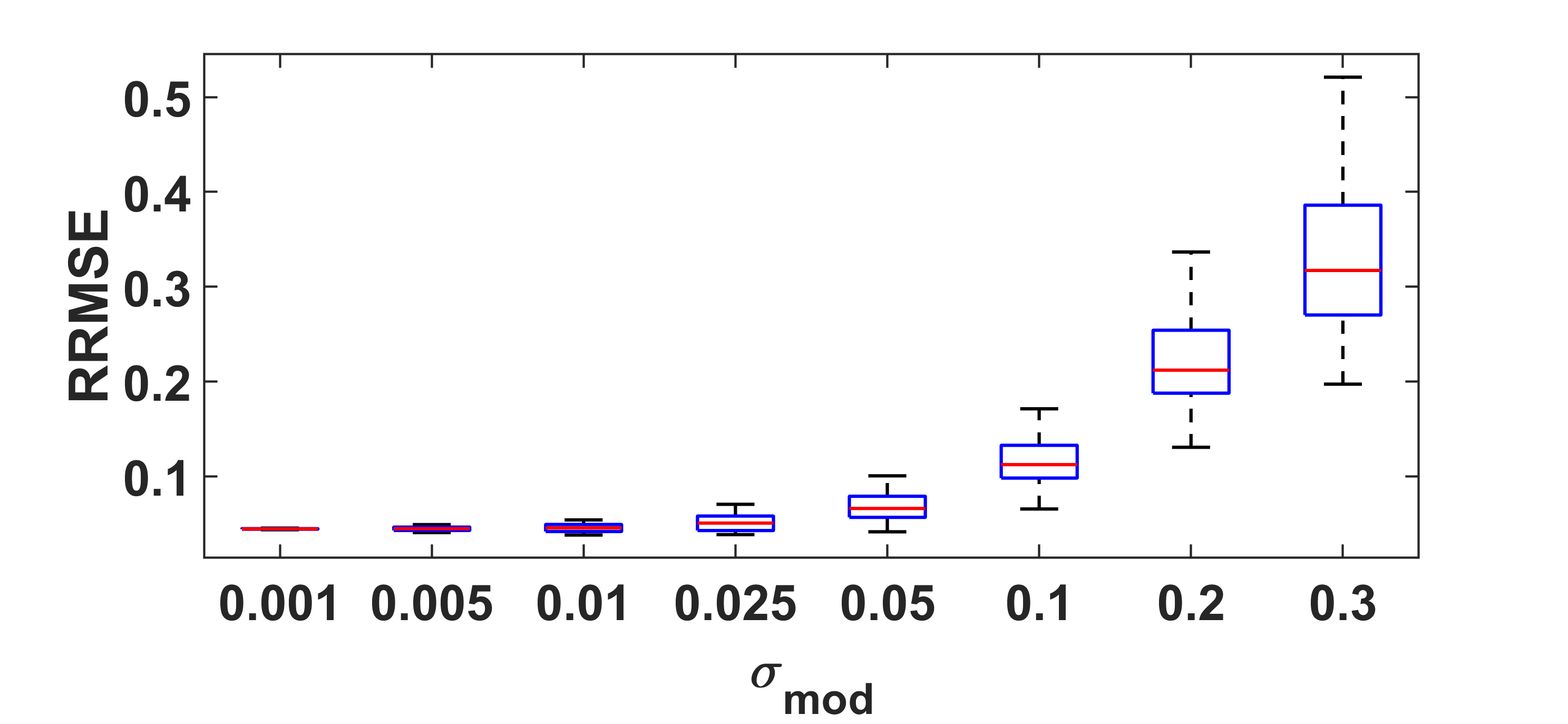}
    \caption{Boxplots for the temperature RRMSE for 50 downscaled fields in the case of only model noise during evaluation. Here, CDAnet was perturbed during training. Here, uncertainties were incorporated to CDAnet's training. During inference, the trained CDAnet model was perturbed then employed to downscale a coarse-scale solution trajectory with $\mathcal{S} = 4$, $\mathcal{T} = 4$, and $Ra=10^5$.}
    \label{fig:sec2_modNoise_1}
\end{figure}

\begin{figure}[!htbp]
    \centering
    \includegraphics[width=\linewidth]{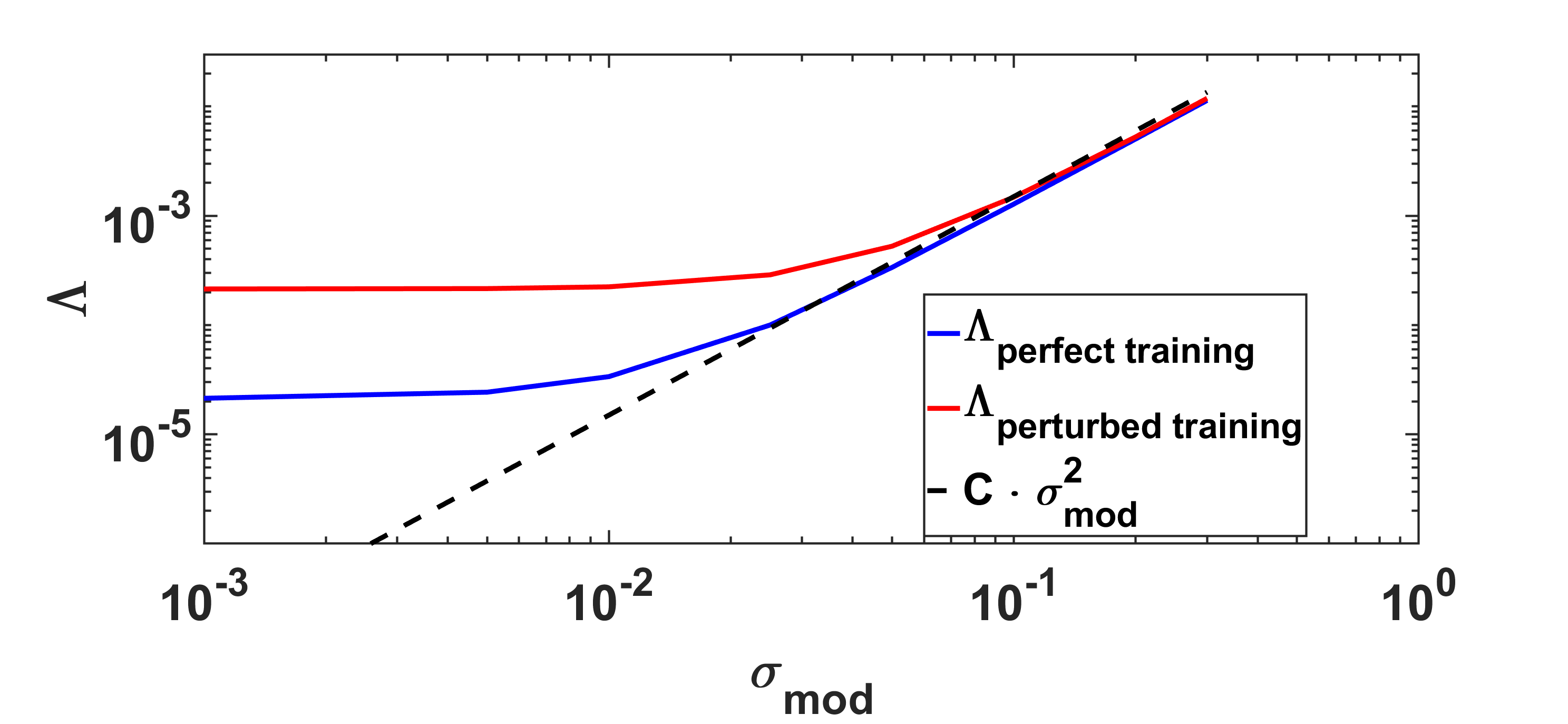}
    \caption{Curves of $\Lambda$ as a function of $\sigma_{mod}$ are shown for the different $\sigma_{obs}$ when CDAnet was trained under perfect and perturbed conditions; $(\mathcal{S}, \mathcal{T}) = (4, 4)$ and $Ra=10^5$.}
    \label{fig:Sec1_2_lambda_sigMod}
\end{figure}


\begin{table}[!htbp]
\centering
\begin{center}
\resizebox{\textwidth}{!}{%
\begin{tabular}{|cc|llllll|}
\hline
\multicolumn{2}{|c|}{\multirow{2}{*}{}} & \multicolumn{6}{c|}{$\sigma_{mod}$} \\ \cline{3-8} 
\multicolumn{2}{|c|}{} & \multicolumn{1}{c|}{0.001} & \multicolumn{1}{c|}{0.005} & \multicolumn{1}{c|}{0.01} & \multicolumn{1}{c|}{0.025} & \multicolumn{1}{c|}{0.05} & \multicolumn{1}{c|}{0.1} \\ \hline
\multicolumn{1}{|c|}{\multirow{8}{*}{$\sigma_{obs}$}} & 0.001 & \multicolumn{1}{l|}{0.044584} & \multicolumn{1}{l|}{0.044728} & \multicolumn{1}{l|}{0.04546} & \multicolumn{1}{l|}{0.051092} & \multicolumn{1}{l|}{0.068658} & 0.115356 \\ \cline{2-8} 
\multicolumn{1}{|c|}{} & 0.0025 & \multicolumn{1}{l|}{0.045103} & \multicolumn{1}{l|}{0.045247} & \multicolumn{1}{l|}{0.045972} & \multicolumn{1}{l|}{0.051555} & \multicolumn{1}{l|}{0.069007} & 0.115563 \\ \cline{2-8} 
\multicolumn{1}{|c|}{} & 0.005 & \multicolumn{1}{l|}{0.046881} & \multicolumn{1}{l|}{0.047022} & \multicolumn{1}{l|}{0.047727} & \multicolumn{1}{l|}{0.053146} & \multicolumn{1}{l|}{0.070211} & 0.116281 \\ \cline{2-8} 
\multicolumn{1}{|c|}{} & 0.0075 & \multicolumn{1}{l|}{0.04966} & \multicolumn{1}{l|}{0.049795} & \multicolumn{1}{l|}{0.050469} & \multicolumn{1}{l|}{0.055651} & \multicolumn{1}{l|}{0.07214} & 0.11745 \\ \cline{2-8} 
\multicolumn{1}{|c|}{} & 0.01 & \multicolumn{1}{l|}{0.053263} & \multicolumn{1}{l|}{0.05339} & \multicolumn{1}{l|}{0.054025} & \multicolumn{1}{l|}{0.058928} & \multicolumn{1}{l|}{0.074718} & 0.119051 \\ \cline{2-8} 
\multicolumn{1}{|c|}{} & 0.025 & \multicolumn{1}{l|}{0.084621} & \multicolumn{1}{l|}{0.08469} & \multicolumn{1}{l|}{0.085097} & \multicolumn{1}{l|}{0.088376} & \multicolumn{1}{l|}{0.099764} & 0.136299 \\ \cline{2-8} 
\multicolumn{1}{|c|}{} & 0.05 & \multicolumn{1}{l|}{0.14893} & \multicolumn{1}{l|}{0.148936} & \multicolumn{1}{l|}{0.14913} & \multicolumn{1}{l|}{0.15095} & \multicolumn{1}{l|}{0.157867} & 0.183394 \\ \cline{2-8} 
\multicolumn{1}{|c|}{} & 0.1 & \multicolumn{1}{l|}{0.283994} & \multicolumn{1}{l|}{0.283913} & \multicolumn{1}{l|}{0.283917} & \multicolumn{1}{l|}{0.284624} & \multicolumn{1}{l|}{0.288066} & 0.302802 \\ \hline
\end{tabular}%
}
\caption*{(a) Mean RRMSE}
\end{center}

\begin{center}
\resizebox{\textwidth}{!}{%
\begin{tabular}{|cc|llllll|}
\hline
\multicolumn{2}{|c|}{\multirow{2}{*}{}} & \multicolumn{6}{c|}{$\sigma_{mod}$} \\ \cline{3-8} 
\multicolumn{2}{|c|}{} & \multicolumn{1}{c|}{0.001} & \multicolumn{1}{c|}{0.005} & \multicolumn{1}{c|}{0.01} & \multicolumn{1}{c|}{0.025} & \multicolumn{1}{c|}{0.05} & \multicolumn{1}{c|}{0.1} \\ \hline
\multicolumn{1}{|c|}{\multirow{8}{*}{$\sigma_{obs}$}} & 0.001 & \multicolumn{1}{l|}{0.000401} & \multicolumn{1}{l|}{0.001991} & \multicolumn{1}{l|}{0.003893} & \multicolumn{1}{l|}{0.008748} & \multicolumn{1}{l|}{0.014615} & 0.024394 \\ \cline{2-8} 
\multicolumn{1}{|c|}{} & 0.0025 & \multicolumn{1}{l|}{0.000396} & \multicolumn{1}{l|}{0.001968} & \multicolumn{1}{l|}{0.003853} & \multicolumn{1}{l|}{0.00868} & \multicolumn{1}{l|}{0.014554} & 0.024363 \\ \cline{2-8} 
\multicolumn{1}{|c|}{} & 0.005 & \multicolumn{1}{l|}{0.000381} & \multicolumn{1}{l|}{0.001898} & \multicolumn{1}{l|}{0.003722} & \multicolumn{1}{l|}{0.008452} & \multicolumn{1}{l|}{0.014346} & 0.024248 \\ \cline{2-8} 
\multicolumn{1}{|c|}{} & 0.0075 & \multicolumn{1}{l|}{0.000359} & \multicolumn{1}{l|}{0.001799} & \multicolumn{1}{l|}{0.003539} & \multicolumn{1}{l|}{0.008122} & \multicolumn{1}{l|}{0.014025} & 0.024058 \\ \cline{2-8} 
\multicolumn{1}{|c|}{} & 0.01 & \multicolumn{1}{l|}{0.000336} & \multicolumn{1}{l|}{0.00169} & \multicolumn{1}{l|}{0.003334} & \multicolumn{1}{l|}{0.007738} & \multicolumn{1}{l|}{0.013623} & 0.023797 \\ \cline{2-8} 
\multicolumn{1}{|c|}{} & 0.025 & \multicolumn{1}{l|}{0.000235} & \multicolumn{1}{l|}{0.001155} & \multicolumn{1}{l|}{0.002304} & \multicolumn{1}{l|}{0.005625} & \multicolumn{1}{l|}{0.010915} & 0.021451 \\ \cline{2-8} 
\multicolumn{1}{|c|}{} & 0.05 & \multicolumn{1}{l|}{0.000217} & \multicolumn{1}{l|}{0.000844} & \multicolumn{1}{l|}{0.001678} & \multicolumn{1}{l|}{0.004159} & \multicolumn{1}{l|}{0.008384} & 0.017897 \\ \cline{2-8} 
\multicolumn{1}{|c|}{} & 0.1 & \multicolumn{1}{l|}{0.000305} & \multicolumn{1}{l|}{0.000854} & \multicolumn{1}{l|}{0.001672} & \multicolumn{1}{l|}{0.004138} & \multicolumn{1}{l|}{0.008216} & 0.016838 \\ \hline
\end{tabular}%
}
\caption*{(b) Standard Deviation of RRMSE}
\end{center}

 \caption{Temperature RRMSE: (a) mean and (b) standard deviation for different combinations of $\sigma_{obs}$ and $\sigma_{mod}$, for the case when observational and model uncertainties were incorporated during training. During inference, noisy CDAnet models were used to downscale an ensemble of noisy representation of a coarse-scale reference trajectory with $\mathcal{S} = 4$, $\mathcal{T} = 4$, and $Ra=10^5$. \label{fig:sec2_aggNoise_meanStd}}
\end{table}



\clearpage
\newpage


\begin{figure}[!htbp]
    \centering
    \includegraphics[width=\linewidth]{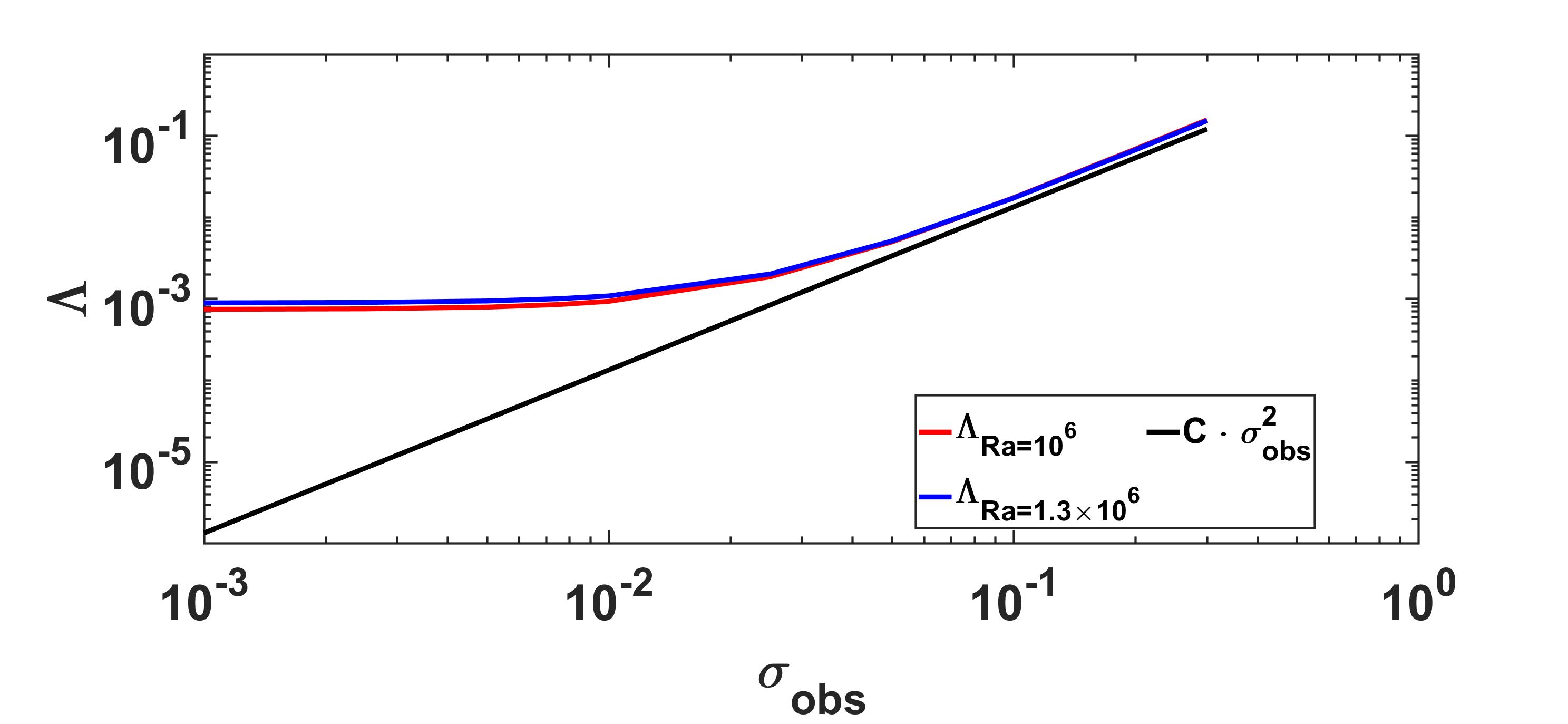}
    \caption{Curves of $\Lambda$ as a function of $\sigma_{obs}$ are shown for both $Ra$ numbers. The results are presented for a CDAnet model trained with imperfectly downscaled fields with $(\mathcal{S}, \mathcal{T}) = (4, 4)$.}
    \label{fig:Sec3_lambda_sigObs}
\end{figure}


\begin{figure}[!htbp]
    \centering
    \subfloat[$Ra=10^6$ - Temperature]{\includegraphics[width=0.49\linewidth]{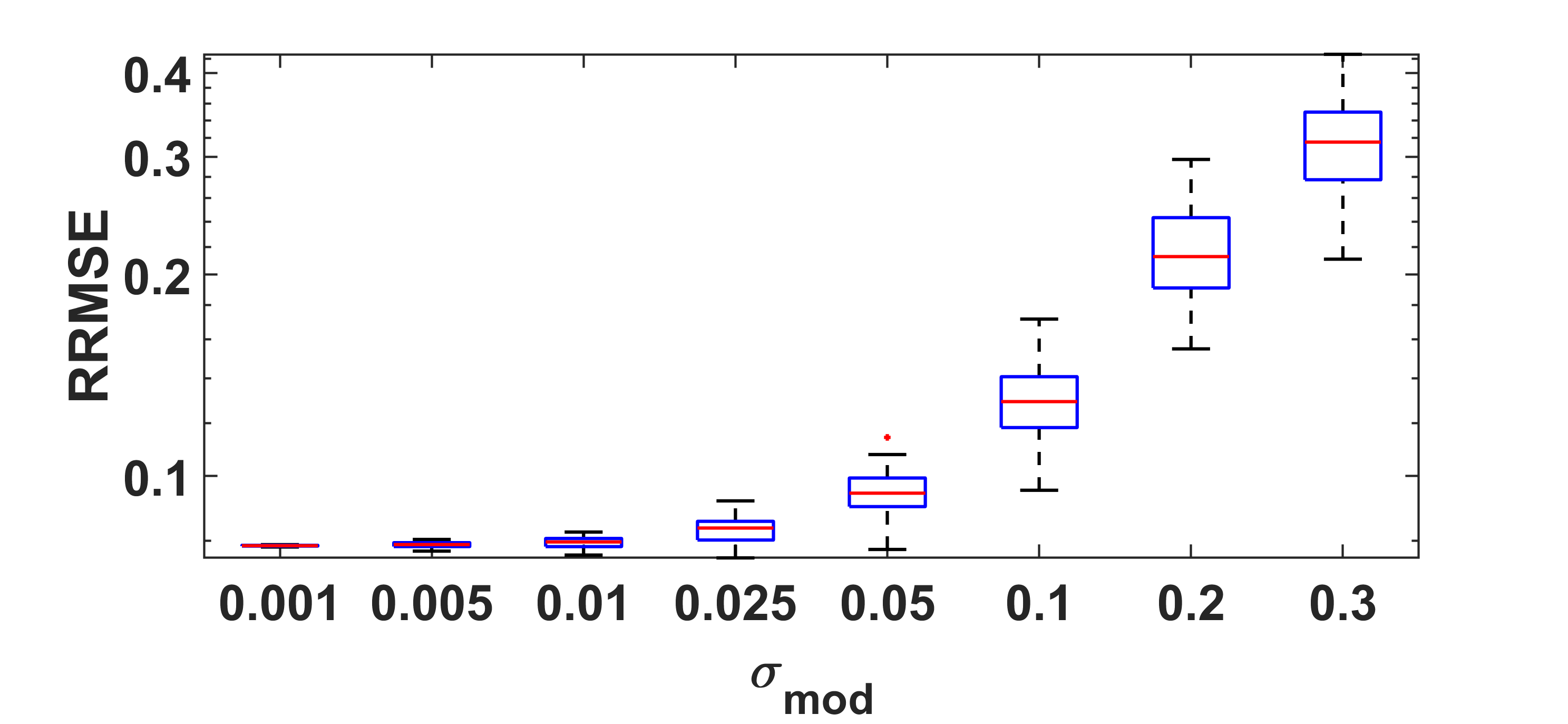}} \,
    \subfloat[$Ra=1.3\times10^6$ - Temperature]{\includegraphics[width=0.49\linewidth]{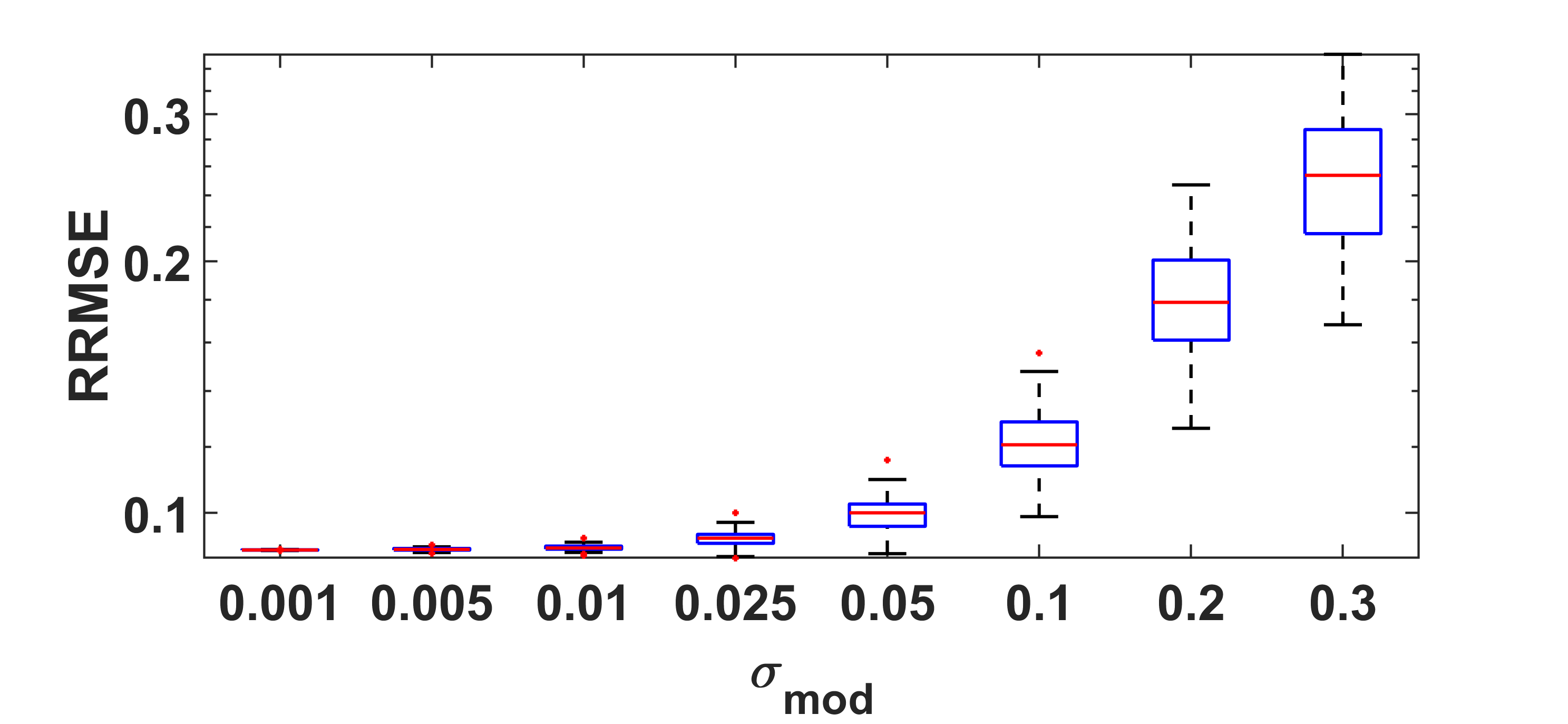}} \,
    \subfloat[$Ra=10^6$ - u-velocity]{\includegraphics[width=0.49\linewidth]{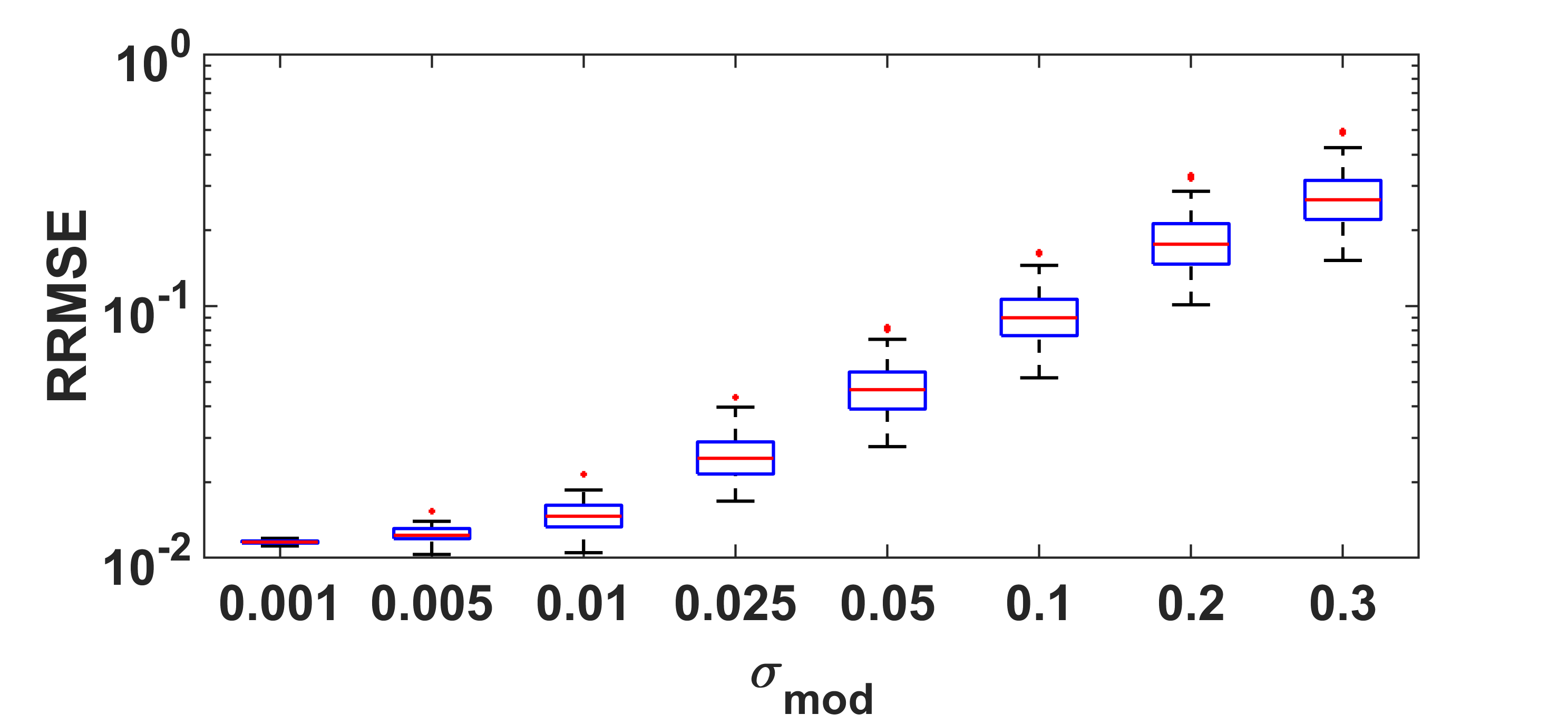}} \,
    \subfloat[$Ra=1.3\times10^6$ - u-velocity]{\includegraphics[width=0.49\linewidth]{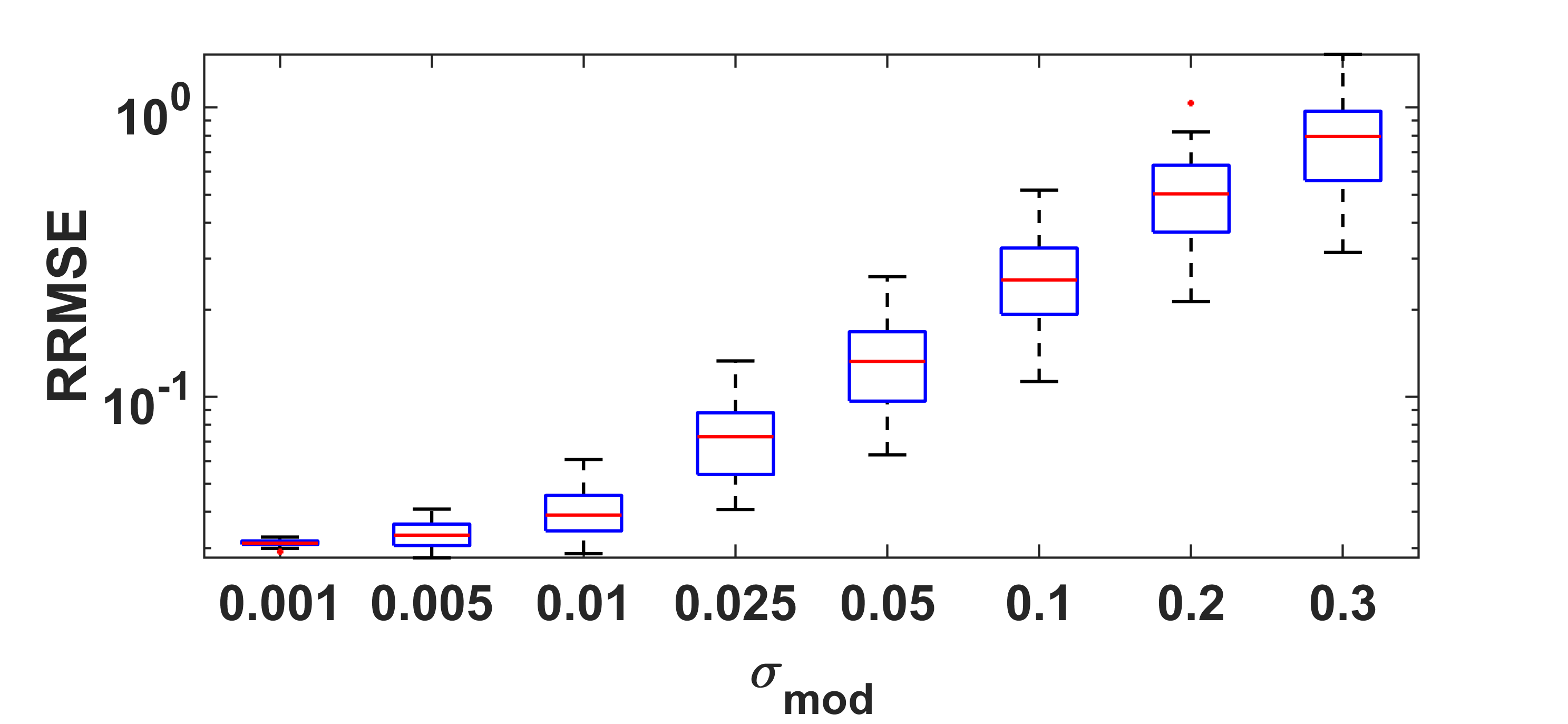}} \,
    \caption{Boxplots of the (a, b) temperature and (c, d) horizontal velocity RRMSE of 50 CDAnet-downscaled realizations with $Ra$ of (a, c) $10^6$ and (b, d) $1.3 \times 10^6$, in the case of only model noise in CDAnet evaluation. In both cases, $\mathcal{S} = 4$ and $\mathcal{T} = 4$, where the CDAnet model was trained with imperfectly downscaled fields.}
    \label{fig:sec3_modNoise_1}
\end{figure}

\begin{figure}[!htbp]
    \centering
    \includegraphics[width=\linewidth]{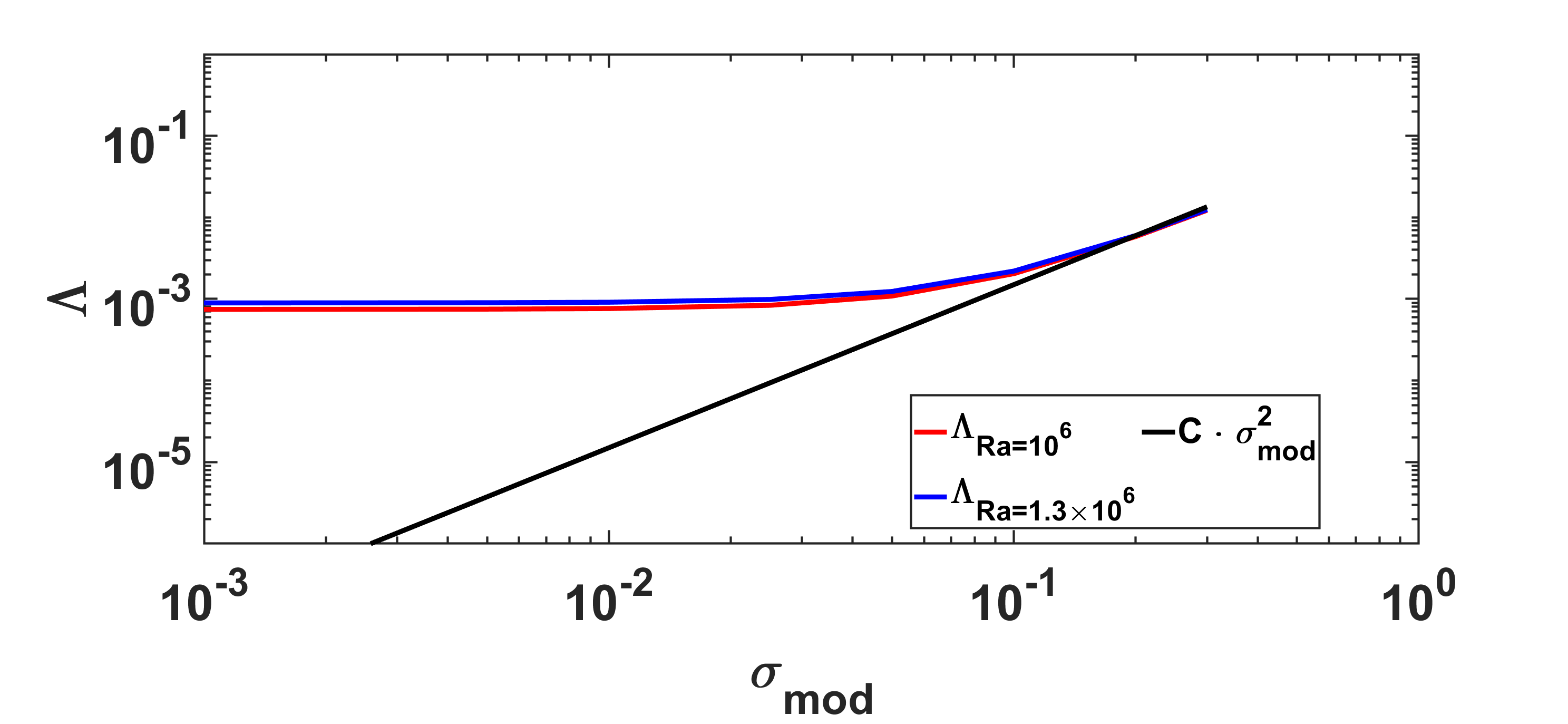}
    \caption{Curves of $\Lambda$ as a function of $\sigma_{mod}$ are shown for both $Ra$ numbers. The results are presented for a CDAnet model trained with imperfectly downscaled fields with $(\mathcal{S}, \mathcal{T}) = (4, 4)$.}
    \label{fig:Sec3_lambda_sigMod}
\end{figure}

\clearpage
\newpage



\begin{table}[!htbp]
\centering
\begin{center}
\resizebox{0.75 \textwidth}{!}{%
\begin{tabular}{|cc|llllll|}
\hline
\multicolumn{2}{|c|}{\multirow{2}{*}{}} & \multicolumn{6}{c|}{$\sigma_{mod}$} \\ \cline{3-8} 
\multicolumn{2}{|c|}{} & \multicolumn{1}{c|}{0.001} & \multicolumn{1}{c|}{0.005} & \multicolumn{1}{c|}{0.01} & \multicolumn{1}{c|}{0.025} & \multicolumn{1}{c|}{0.05} & \multicolumn{1}{c|}{0.1} \\ \hline
\multicolumn{1}{|c|}{\multirow{8}{*}{$\sigma_{obs}$}} & 0.001 & \multicolumn{1}{l|}{0.078814} & \multicolumn{1}{l|}{0.079076} & \multicolumn{1}{l|}{0.079697} & \multicolumn{1}{l|}{0.083418} & \multicolumn{1}{l|}{0.094819} & 0.12957 \\ \cline{2-8} 
\multicolumn{1}{|c|}{} & 0.0025 & \multicolumn{1}{l|}{0.079365} & \multicolumn{1}{l|}{0.079627} & \multicolumn{1}{l|}{0.080244} & \multicolumn{1}{l|}{0.083942} & \multicolumn{1}{l|}{0.095283} & 0.129919 \\ \cline{2-8} 
\multicolumn{1}{|c|}{} & 0.005 & \multicolumn{1}{l|}{0.081283} & \multicolumn{1}{l|}{0.081539} & \multicolumn{1}{l|}{0.082143} & \multicolumn{1}{l|}{0.085761} & \multicolumn{1}{l|}{0.096894} & 0.13112 \\ \cline{2-8} 
\multicolumn{1}{|c|}{} & 0.0075 & \multicolumn{1}{l|}{0.084308} & \multicolumn{1}{l|}{0.084555} & \multicolumn{1}{l|}{0.085138} & \multicolumn{1}{l|}{0.088633} & \multicolumn{1}{l|}{0.099444} & 0.133025 \\ \cline{2-8} 
\multicolumn{1}{|c|}{} & 0.01 & \multicolumn{1}{l|}{0.08828} & \multicolumn{1}{l|}{0.088517} & \multicolumn{1}{l|}{0.089074} & \multicolumn{1}{l|}{0.092416} & \multicolumn{1}{l|}{0.102823} & 0.135575 \\ \cline{2-8} 
\multicolumn{1}{|c|}{} & 0.025 & \multicolumn{1}{l|}{0.124829} & \multicolumn{1}{l|}{0.124996} & \multicolumn{1}{l|}{0.125386} & \multicolumn{1}{l|}{0.127745} & \multicolumn{1}{l|}{0.135389} & 0.161621 \\ \cline{2-8} 
\multicolumn{1}{|c|}{} & 0.05 & \multicolumn{1}{l|}{0.204696} & \multicolumn{1}{l|}{0.204796} & \multicolumn{1}{l|}{0.20503} & \multicolumn{1}{l|}{0.206458} & \multicolumn{1}{l|}{0.211212} & 0.228884 \\ \cline{2-8} 
\multicolumn{1}{|c|}{} & 0.1 & \multicolumn{1}{l|}{0.382146} & \multicolumn{1}{l|}{0.382173} & \multicolumn{1}{l|}{0.382271} & \multicolumn{1}{l|}{0.382991} & \multicolumn{1}{l|}{0.385588} & 0.39587 \\ \hline
\end{tabular}%
}
\caption*{(a) Mean RRMSE - $Ra = 10^6$}
\end{center}

\begin{center}
\resizebox{0.75 \textwidth}{!}{%
\begin{tabular}{|cc|llllll|}
\hline
\multicolumn{2}{|c|}{\multirow{2}{*}{}} & \multicolumn{6}{c|}{$\sigma_{mod}$} \\ \cline{3-8} 
\multicolumn{2}{|c|}{} & \multicolumn{1}{c|}{0.001} & \multicolumn{1}{c|}{0.005} & \multicolumn{1}{c|}{0.01} & \multicolumn{1}{c|}{0.025} & \multicolumn{1}{c|}{0.05} & \multicolumn{1}{c|}{0.1} \\ \hline
\multicolumn{1}{|c|}{\multirow{8}{*}{$\sigma_{obs}$}} & 0.001 & \multicolumn{1}{l|}{0.000162} & \multicolumn{1}{l|}{0.000801} & \multicolumn{1}{l|}{0.001582} & \multicolumn{1}{l|}{0.003863} & \multicolumn{1}{l|}{0.007741} & 0.016528 \\ \cline{2-8} 
\multicolumn{1}{|c|}{} & 0.0025 & \multicolumn{1}{l|}{0.000164} & \multicolumn{1}{l|}{0.000798} & \multicolumn{1}{l|}{0.001575} & \multicolumn{1}{l|}{0.003842} & \multicolumn{1}{l|}{0.007705} & 0.016479 \\ \cline{2-8} 
\multicolumn{1}{|c|}{} & 0.005 & \multicolumn{1}{l|}{0.000165} & \multicolumn{1}{l|}{0.000783} & \multicolumn{1}{l|}{0.001543} & \multicolumn{1}{l|}{0.003768} & \multicolumn{1}{l|}{0.007585} & 0.016318 \\ \cline{2-8} 
\multicolumn{1}{|c|}{} & 0.0075 & \multicolumn{1}{l|}{0.000164} & \multicolumn{1}{l|}{0.00076} & \multicolumn{1}{l|}{0.001495} & \multicolumn{1}{l|}{0.003658} & \multicolumn{1}{l|}{0.007404} & 0.016072 \\ \cline{2-8} 
\multicolumn{1}{|c|}{} & 0.01 & \multicolumn{1}{l|}{0.000165} & \multicolumn{1}{l|}{0.000733} & \multicolumn{1}{l|}{0.001439} & \multicolumn{1}{l|}{0.003524} & \multicolumn{1}{l|}{0.007177} & 0.015755 \\ \cline{2-8} 
\multicolumn{1}{|c|}{} & 0.025 & \multicolumn{1}{l|}{0.000174} & \multicolumn{1}{l|}{0.00057} & \multicolumn{1}{l|}{0.00109} & \multicolumn{1}{l|}{0.00265} & \multicolumn{1}{l|}{0.005523} & 0.013065 \\ \cline{2-8} 
\multicolumn{1}{|c|}{} & 0.05 & \multicolumn{1}{l|}{0.000204} & \multicolumn{1}{l|}{0.000425} & \multicolumn{1}{l|}{0.000761} & \multicolumn{1}{l|}{0.001779} & \multicolumn{1}{l|}{0.003621} & 0.008945 \\ \cline{2-8} 
\multicolumn{1}{|c|}{} & 0.1 & \multicolumn{1}{l|}{0.000287} & \multicolumn{1}{l|}{0.000507} & \multicolumn{1}{l|}{0.000888} & \multicolumn{1}{l|}{0.002077} & \multicolumn{1}{l|}{0.004027} & 0.008322 \\ \hline
\end{tabular}%
}
\caption*{(b) Standard Deviation of RRMSE - $Ra = 10^6$}
\end{center}

\begin{center}
\resizebox{0.75 \textwidth}{!}{%
\begin{tabular}{|cc|llllll|}
\hline
\multicolumn{2}{|c|}{\multirow{2}{*}{}} & \multicolumn{6}{c|}{$\sigma_{mod}$} \\ \cline{3-8} 
\multicolumn{2}{|c|}{} & \multicolumn{1}{c|}{0.001} & \multicolumn{1}{c|}{0.005} & \multicolumn{1}{c|}{0.01} & \multicolumn{1}{c|}{0.025} & \multicolumn{1}{c|}{0.05} & \multicolumn{1}{c|}{0.1} \\ \hline
\multicolumn{1}{|c|}{\multirow{8}{*}{$\sigma_{obs}$}} & 0.001 & \multicolumn{1}{l|}{0.090367} & \multicolumn{1}{l|}{0.090528} & \multicolumn{1}{l|}{0.090915} & \multicolumn{1}{l|}{0.093212} & \multicolumn{1}{l|}{0.100169} & 0.121875 \\ \cline{2-8} 
\multicolumn{1}{|c|}{} & 0.0025 & \multicolumn{1}{l|}{0.090832} & \multicolumn{1}{l|}{0.090992} & \multicolumn{1}{l|}{0.091375} & \multicolumn{1}{l|}{0.093651} & \multicolumn{1}{l|}{0.100561} & 0.122181 \\ \cline{2-8} 
\multicolumn{1}{|c|}{} & 0.005 & \multicolumn{1}{l|}{0.092427} & \multicolumn{1}{l|}{0.092581} & \multicolumn{1}{l|}{0.092952} & \multicolumn{1}{l|}{0.095158} & \multicolumn{1}{l|}{0.101908} & 0.123238 \\ \cline{2-8} 
\multicolumn{1}{|c|}{} & 0.0075 & \multicolumn{1}{l|}{0.094903} & \multicolumn{1}{l|}{0.09505} & \multicolumn{1}{l|}{0.095403} & \multicolumn{1}{l|}{0.097509} & \multicolumn{1}{l|}{0.104027} & 0.124918 \\ \cline{2-8} 
\multicolumn{1}{|c|}{} & 0.01 & \multicolumn{1}{l|}{0.098109} & \multicolumn{1}{l|}{0.098249} & \multicolumn{1}{l|}{0.098581} & \multicolumn{1}{l|}{0.100576} & \multicolumn{1}{l|}{0.106818} & 0.127164 \\ \cline{2-8} 
\multicolumn{1}{|c|}{} & 0.025 & \multicolumn{1}{l|}{0.126784} & \multicolumn{1}{l|}{0.126883} & \multicolumn{1}{l|}{0.127115} & \multicolumn{1}{l|}{0.128533} & \multicolumn{1}{l|}{0.133171} & 0.149638 \\ \cline{2-8} 
\multicolumn{1}{|c|}{} & 0.05 & \multicolumn{1}{l|}{0.190569} & \multicolumn{1}{l|}{0.190626} & \multicolumn{1}{l|}{0.190769} & \multicolumn{1}{l|}{0.19167} & \multicolumn{1}{l|}{0.194733} & 0.206342 \\ \cline{2-8} 
\multicolumn{1}{|c|}{} & 0.1 & \multicolumn{1}{l|}{0.335152} & \multicolumn{1}{l|}{0.335159} & \multicolumn{1}{l|}{0.335213} & \multicolumn{1}{l|}{0.335682} & \multicolumn{1}{l|}{0.337489} & 0.344922 \\ \hline
\end{tabular}%
}
\caption*{(c) Mean RRMSE - $Ra = 1.3 \times 10^6$}
\end{center}

\begin{center}
\resizebox{0.75 \textwidth}{!}{%
\begin{tabular}{|cc|llllll|}
\hline
\multicolumn{2}{|c|}{\multirow{2}{*}{}} & \multicolumn{6}{c|}{$\sigma_{mod}$} \\ \cline{3-8} 
\multicolumn{2}{|c|}{} & \multicolumn{1}{c|}{0.001} & \multicolumn{1}{c|}{0.005} & \multicolumn{1}{c|}{0.01} & \multicolumn{1}{c|}{0.025} & \multicolumn{1}{c|}{0.05} & \multicolumn{1}{c|}{0.1} \\ \hline
\multicolumn{1}{|c|}{\multirow{8}{*}{$\sigma_{obs}$}} & 0.001 & \multicolumn{1}{l|}{7.87E-05} & \multicolumn{1}{l|}{0.00039} & \multicolumn{1}{l|}{0.00079} & \multicolumn{1}{l|}{0.002164} & \multicolumn{1}{l|}{0.005069} & 0.012335 \\ \cline{2-8} 
\multicolumn{1}{|c|}{} & 0.0025 & \multicolumn{1}{l|}{8.16E-05} & \multicolumn{1}{l|}{0.000391} & \multicolumn{1}{l|}{0.000789} & \multicolumn{1}{l|}{0.002151} & \multicolumn{1}{l|}{0.005041} & 0.012294 \\ \cline{2-8} 
\multicolumn{1}{|c|}{} & 0.005 & \multicolumn{1}{l|}{8.91E-05} & \multicolumn{1}{l|}{0.000392} & \multicolumn{1}{l|}{0.000785} & \multicolumn{1}{l|}{0.002114} & \multicolumn{1}{l|}{0.00495} & 0.012155 \\ \cline{2-8} 
\multicolumn{1}{|c|}{} & 0.0075 & \multicolumn{1}{l|}{9.89E-05} & \multicolumn{1}{l|}{0.000392} & \multicolumn{1}{l|}{0.000777} & \multicolumn{1}{l|}{0.002065} & \multicolumn{1}{l|}{0.00482} & 0.011946 \\ \cline{2-8} 
\multicolumn{1}{|c|}{} & 0.01 & \multicolumn{1}{l|}{0.00011} & \multicolumn{1}{l|}{0.000391} & \multicolumn{1}{l|}{0.000766} & \multicolumn{1}{l|}{0.00201} & \multicolumn{1}{l|}{0.004667} & 0.011678 \\ \cline{2-8} 
\multicolumn{1}{|c|}{} & 0.025 & \multicolumn{1}{l|}{0.000164} & \multicolumn{1}{l|}{0.000377} & \multicolumn{1}{l|}{0.000691} & \multicolumn{1}{l|}{0.00169} & \multicolumn{1}{l|}{0.003721} & 0.009598 \\ \cline{2-8} 
\multicolumn{1}{|c|}{} & 0.05 & \multicolumn{1}{l|}{0.000237} & \multicolumn{1}{l|}{0.000383} & \multicolumn{1}{l|}{0.000626} & \multicolumn{1}{l|}{0.0014} & \multicolumn{1}{l|}{0.002814} & 0.006863 \\ \cline{2-8} 
\multicolumn{1}{|c|}{} & 0.1 & \multicolumn{1}{l|}{0.000339} & \multicolumn{1}{l|}{0.00048} & \multicolumn{1}{l|}{0.000767} & \multicolumn{1}{l|}{0.001738} & \multicolumn{1}{l|}{0.003377} & 0.006899 \\ \hline
\end{tabular}%
}
\caption*{(d) Standard Deviation of RRMSE - $Ra = 1.3\times 10^6$}
\end{center}

\caption{Mean and standard deviation of the temperature RRMSE of 50 CDAnet-downscaled realizations, estimated with $Ra=10^6$ and $1.3\times10^6$, for the case when CDAnet was trained with imperfectly downscaled fields. During inference, noisy CDAnet models were used to downscale an ensemble of noisy representation of a coarse-scale reference trajectory with $\mathcal{S} = 4$ and $\mathcal{T} = 4$. 
\label{fig:sec3_combinedNoise}}

\end{table}



\begin{table}[!htbp]
\centering
\begin{center}
\resizebox{0.75 \textwidth}{!}{%
\begin{tabular}{|cc|llllll|}
\hline
\multicolumn{2}{|c|}{\multirow{2}{*}{}} & \multicolumn{6}{c|}{$\sigma_{mod}$} \\ \cline{3-8} 
\multicolumn{2}{|c|}{} & \multicolumn{1}{c|}{0.001} & \multicolumn{1}{c|}{0.005} & \multicolumn{1}{c|}{0.01} & \multicolumn{1}{c|}{0.025} & \multicolumn{1}{c|}{0.05} & \multicolumn{1}{c|}{0.1} \\ \hline
\multicolumn{1}{|c|}{\multirow{8}{*}{$\sigma_{obs}$}} & 0.001 & \multicolumn{1}{l|}{0.011919456} & \multicolumn{1}{l|}{0.012755389} & \multicolumn{1}{l|}{0.015099953} & \multicolumn{1}{l|}{0.026116933} & \multicolumn{1}{l|}{0.047902629} & 0.093287086 \\ \cline{2-8} 
\multicolumn{1}{|c|}{} & 0.0025 & \multicolumn{1}{l|}{0.013710172} & \multicolumn{1}{l|}{0.014447563} & \multicolumn{1}{l|}{0.016566833} & \multicolumn{1}{l|}{0.027012605} & \multicolumn{1}{l|}{0.048402335} & 0.09354413 \\ \cline{2-8} 
\multicolumn{1}{|c|}{} & 0.005 & \multicolumn{1}{l|}{0.018732265} & \multicolumn{1}{l|}{0.019283435} & \multicolumn{1}{l|}{0.020939422} & \multicolumn{1}{l|}{0.029958939} & \multicolumn{1}{l|}{0.050135397} & 0.094453949 \\ \cline{2-8} 
\multicolumn{1}{|c|}{} & 0.0075 & \multicolumn{1}{l|}{0.02493467} & \multicolumn{1}{l|}{0.025350852} & \multicolumn{1}{l|}{0.026641814} & \multicolumn{1}{l|}{0.034254367} & \multicolumn{1}{l|}{0.052871171} & 0.095946121 \\ \cline{2-8} 
\multicolumn{1}{|c|}{} & 0.01 & \multicolumn{1}{l|}{0.031633286} & \multicolumn{1}{l|}{0.031959217} & \multicolumn{1}{l|}{0.032994061} & \multicolumn{1}{l|}{0.039443379} & \multicolumn{1}{l|}{0.056442934} & 0.097988169 \\ \cline{2-8} 
\multicolumn{1}{|c|}{} & 0.025 & \multicolumn{1}{l|}{0.074480419} & \multicolumn{1}{l|}{0.074597852} & \multicolumn{1}{l|}{0.075023896} & \multicolumn{1}{l|}{0.078073881} & \multicolumn{1}{l|}{0.088094944} & 0.119418037 \\ \cline{2-8} 
\multicolumn{1}{|c|}{} & 0.05 & \multicolumn{1}{l|}{0.147611126} & \multicolumn{1}{l|}{0.147634009} & \multicolumn{1}{l|}{0.147804834} & \multicolumn{1}{l|}{0.149256385} & \multicolumn{1}{l|}{0.154646798} & 0.174664434 \\ \cline{2-8} 
\multicolumn{1}{|c|}{} & 0.1 & \multicolumn{1}{l|}{0.294361047} & \multicolumn{1}{l|}{0.294302245} & \multicolumn{1}{l|}{0.294302361} & \multicolumn{1}{l|}{0.294796246} & \multicolumn{1}{l|}{0.297247251} & 0.307952166 \\ \hline
\end{tabular}%
}
\caption*{(a) Mean RRMSE - $Ra = 10^6$}
\end{center}

\begin{center}
\resizebox{0.75 \textwidth}{!}{%
\begin{tabular}{|cc|llllll|}
\hline
\multicolumn{2}{|c|}{\multirow{2}{*}{}} & \multicolumn{6}{c|}{$\sigma_{mod}$} \\ \cline{3-8} 
\multicolumn{2}{|c|}{} & \multicolumn{1}{c|}{0.001} & \multicolumn{1}{c|}{0.005} & \multicolumn{1}{c|}{0.01} & \multicolumn{1}{c|}{0.025} & \multicolumn{1}{c|}{0.05} & \multicolumn{1}{c|}{0.1} \\ \hline
\multicolumn{1}{|c|}{\multirow{8}{*}{$\sigma_{obs}$}} & 0.001 & \multicolumn{1}{l|}{0.000176886} & \multicolumn{1}{l|}{0.000903422} & \multicolumn{1}{l|}{0.001978722} & \multicolumn{1}{l|}{0.005687393} & \multicolumn{1}{l|}{0.012022838} & 0.024743948 \\ \cline{2-8} 
\multicolumn{1}{|c|}{} & 0.0025 & \multicolumn{1}{l|}{0.00015475} & \multicolumn{1}{l|}{0.000801805} & \multicolumn{1}{l|}{0.00181445} & \multicolumn{1}{l|}{0.005518198} & \multicolumn{1}{l|}{0.01190941} & 0.024679544 \\ \cline{2-8} 
\multicolumn{1}{|c|}{} & 0.005 & \multicolumn{1}{l|}{0.000116081} & \multicolumn{1}{l|}{0.000611071} & \multicolumn{1}{l|}{0.00146083} & \multicolumn{1}{l|}{0.005033073} & \multicolumn{1}{l|}{0.011535829} & 0.024456279 \\ \cline{2-8} 
\multicolumn{1}{|c|}{} & 0.0075 & \multicolumn{1}{l|}{9.12E-05} & \multicolumn{1}{l|}{0.00047541} & \multicolumn{1}{l|}{0.001169563} & \multicolumn{1}{l|}{0.004467997} & \multicolumn{1}{l|}{0.01099874} & 0.024102467 \\ \cline{2-8} 
\multicolumn{1}{|c|}{} & 0.01 & \multicolumn{1}{l|}{7.72E-05} & \multicolumn{1}{l|}{0.000388052} & \multicolumn{1}{l|}{0.000961383} & \multicolumn{1}{l|}{0.003936638} & \multicolumn{1}{l|}{0.010375775} & 0.023639241 \\ \cline{2-8} 
\multicolumn{1}{|c|}{} & 0.025 & \multicolumn{1}{l|}{7.52E-05} & \multicolumn{1}{l|}{0.000238038} & \multicolumn{1}{l|}{0.000516859} & \multicolumn{1}{l|}{0.002114668} & \multicolumn{1}{l|}{0.006940451} & 0.019777171 \\ \cline{2-8} 
\multicolumn{1}{|c|}{} & 0.05 & \multicolumn{1}{l|}{0.000126324} & \multicolumn{1}{l|}{0.000289137} & \multicolumn{1}{l|}{0.000539048} & \multicolumn{1}{l|}{0.001493001} & \multicolumn{1}{l|}{0.00432522} & 0.014150733 \\ \cline{2-8} 
\multicolumn{1}{|c|}{} & 0.1 & \multicolumn{1}{l|}{0.000238218} & \multicolumn{1}{l|}{0.000497849} & \multicolumn{1}{l|}{0.00090908} & \multicolumn{1}{l|}{0.002165959} & \multicolumn{1}{l|}{0.004396332} & 0.010594622 \\ \hline
\end{tabular}%
}
\caption*{(b) Standard Deviation of RRMSE - $Ra = 10^6$}
\end{center}

\begin{center}
\resizebox{0.75 \textwidth}{!}{%
\begin{tabular}{|cc|llllll|}
\hline
\multicolumn{2}{|c|}{\multirow{2}{*}{}} & \multicolumn{6}{c|}{$\sigma_{mod}$} \\ \cline{3-8} 
\multicolumn{2}{|c|}{} & \multicolumn{1}{c|}{0.001} & \multicolumn{1}{c|}{0.005} & \multicolumn{1}{c|}{0.01} & \multicolumn{1}{c|}{0.025} & \multicolumn{1}{c|}{0.05} & \multicolumn{1}{c|}{0.1} \\ \hline
\multicolumn{1}{|c|}{\multirow{8}{*}{$\sigma_{obs}$}} & 0.001 & \multicolumn{1}{l|}{0.034231} & \multicolumn{1}{l|}{0.036408} & \multicolumn{1}{l|}{0.042863} & \multicolumn{1}{l|}{0.073424} & \multicolumn{1}{l|}{0.134124} & 0.260909 \\ \cline{2-8} 
\multicolumn{1}{|c|}{} & 0.0025 & \multicolumn{1}{l|}{0.046554} & \multicolumn{1}{l|}{0.04825} & \multicolumn{1}{l|}{0.05344} & \multicolumn{1}{l|}{0.080456} & \multicolumn{1}{l|}{0.138306} & 0.263152 \\ \cline{2-8} 
\multicolumn{1}{|c|}{} & 0.005 & \multicolumn{1}{l|}{0.075426} & \multicolumn{1}{l|}{0.076542} & \multicolumn{1}{l|}{0.080069} & \multicolumn{1}{l|}{0.100896} & \multicolumn{1}{l|}{0.151908} & 0.270954 \\ \cline{2-8} 
\multicolumn{1}{|c|}{} & 0.0075 & \multicolumn{1}{l|}{0.107488} & \multicolumn{1}{l|}{0.108292} & \multicolumn{1}{l|}{0.110871} & \multicolumn{1}{l|}{0.127189} & \multicolumn{1}{l|}{0.171406} & 0.283103 \\ \cline{2-8} 
\multicolumn{1}{|c|}{} & 0.01 & \multicolumn{1}{l|}{0.140477} & \multicolumn{1}{l|}{0.141098} & \multicolumn{1}{l|}{0.143105} & \multicolumn{1}{l|}{0.156315} & \multicolumn{1}{l|}{0.194724} & 0.298742 \\ \cline{2-8} 
\multicolumn{1}{|c|}{} & 0.025 & \multicolumn{1}{l|}{0.3418} & \multicolumn{1}{l|}{0.342023} & \multicolumn{1}{l|}{0.342824} & \multicolumn{1}{l|}{0.348613} & \multicolumn{1}{l|}{0.368447} & 0.436235 \\ \cline{2-8} 
\multicolumn{1}{|c|}{} & 0.05 & \multicolumn{1}{l|}{0.679888} & \multicolumn{1}{l|}{0.679897} & \multicolumn{1}{l|}{0.680179} & \multicolumn{1}{l|}{0.682806} & \multicolumn{1}{l|}{0.692938} & 0.732552 \\ \cline{2-8} 
\multicolumn{1}{|c|}{} & 0.1 & \multicolumn{1}{l|}{1.362806} & \multicolumn{1}{l|}{1.362545} & \multicolumn{1}{l|}{1.362369} & \multicolumn{1}{l|}{1.362833} & \multicolumn{1}{l|}{1.366892} & 1.386924 \\ \hline
\end{tabular}%
}
\caption*{(c) Mean RRMSE - $Ra = 1.3 \times 10^6$}
\end{center}

\begin{center}
\resizebox{0.75 \textwidth}{!}{%
\begin{tabular}{|cc|llllll|}
\hline
\multicolumn{2}{|c|}{\multirow{2}{*}{}} & \multicolumn{6}{c|}{$\sigma_{mod}$} \\ \cline{3-8} 
\multicolumn{2}{|c|}{} & \multicolumn{1}{c|}{0.001} & \multicolumn{1}{c|}{0.005} & \multicolumn{1}{c|}{0.01} & \multicolumn{1}{c|}{0.025} & \multicolumn{1}{c|}{0.05} & \multicolumn{1}{c|}{0.1} \\ \hline
\multicolumn{1}{|c|}{\multirow{8}{*}{$\sigma_{obs}$}} & 0.001 & \multicolumn{1}{l|}{0.000626} & \multicolumn{1}{l|}{0.003122} & \multicolumn{1}{l|}{0.00675} & \multicolumn{1}{l|}{0.020061} & \multicolumn{1}{l|}{0.042568} & 0.086786 \\ \cline{2-8} 
\multicolumn{1}{|c|}{} & 0.0025 & \multicolumn{1}{l|}{0.000436} & \multicolumn{1}{l|}{0.002284} & \multicolumn{1}{l|}{0.005386} & \multicolumn{1}{l|}{0.018278} & \multicolumn{1}{l|}{0.041205} & 0.085972 \\ \cline{2-8} 
\multicolumn{1}{|c|}{} & 0.005 & \multicolumn{1}{l|}{0.000259} & \multicolumn{1}{l|}{0.001386} & \multicolumn{1}{l|}{0.003591} & \multicolumn{1}{l|}{0.014718} & \multicolumn{1}{l|}{0.037543} & 0.083348 \\ \cline{2-8} 
\multicolumn{1}{|c|}{} & 0.0075 & \multicolumn{1}{l|}{0.000202} & \multicolumn{1}{l|}{0.000956} & \multicolumn{1}{l|}{0.002586} & \multicolumn{1}{l|}{0.011849} & \multicolumn{1}{l|}{0.033507} & 0.079737 \\ \cline{2-8} 
\multicolumn{1}{|c|}{} & 0.01 & \multicolumn{1}{l|}{0.000203} & \multicolumn{1}{l|}{0.00073} & \multicolumn{1}{l|}{0.002} & \multicolumn{1}{l|}{0.009761} & \multicolumn{1}{l|}{0.0298} & 0.075699 \\ \cline{2-8} 
\multicolumn{1}{|c|}{} & 0.025 & \multicolumn{1}{l|}{0.000452} & \multicolumn{1}{l|}{0.000638} & \multicolumn{1}{l|}{0.001191} & \multicolumn{1}{l|}{0.004828} & \multicolumn{1}{l|}{0.016804} & 0.05395 \\ \cline{2-8} 
\multicolumn{1}{|c|}{} & 0.05 & \multicolumn{1}{l|}{0.000924} & \multicolumn{1}{l|}{0.001268} & \multicolumn{1}{l|}{0.002011} & \multicolumn{1}{l|}{0.004895} & \multicolumn{1}{l|}{0.012027} & 0.036466 \\ \cline{2-8} 
\multicolumn{1}{|c|}{} & 0.1 & \multicolumn{1}{l|}{0.001888} & \multicolumn{1}{l|}{0.002656} & \multicolumn{1}{l|}{0.004198} & \multicolumn{1}{l|}{0.009469} & \multicolumn{1}{l|}{0.018642} & 0.039134 \\ \hline
\end{tabular}%
}
\caption*{(d) Standard Deviation of RRMSE - $Ra = 1.3\times 10^6$}
\end{center}

\caption{Mean and standard deviation of the u-velocity RRMSE of 50 CDAnet-downscaled realizations, estimated with $Ra=10^6$ and $1.3\times10^6$, for the case when CDAnet was trained with imperfectly downscaled fields. During inference, noisy CDAnet models were used to downscale an ensemble of noisy representation of a coarse-scale reference trajectory with $\mathcal{S} = 4$ and $\mathcal{T} = 4$. \label{fig:sec3_combinedNoiseb}}

\end{table}

\clearpage
\newpage

\appendix

\section{Evaluation Metrics}
\label{ssec:evalMetrics}

The performance of the CDAnet under observational and model uncertainties was assessed using different error metrics. 
The downscaled fields were compared with the reference solution obtained from the noise-free simulations performed on fine mesh. 
The mean absolute error (MAE) between the downscaled and reference solutions was computed as:

\begin{equation}
    AE_f(t) = || f(\bm{x}; t) - \hat{f}(\bm{x}; t) ||_1,
\end{equation}
where $f$ and $\hat{f}$ are generic components of the reference and downscaled solutions, respectively.
Similarly, the root mean squared error (RMSE) is computed as:

\begin{equation}
    RMSE_f(t) = || f(\bm{x}; t) - \hat{f}(\bm{x}; t) ||_2.
\end{equation}
The normalized $\ell_2$ error (RRMSE) was used to select the best trained model \citep{HammoudJAMES2022} and is computed as:

\begin{equation}
    RRMSE_f(t) = \frac{|| f(\bm{x}; t) - \hat{f}(\bm{x}; t) ||_2}{|| f(\bm{x}; t) ||_2}.
\end{equation}

When an ensemble of downscaled solutions is considered, the above-mentioned metrics were computed for individual members of the ensemble, and the statistics of the resulting distribution are consequently analyzed. 
In particular, as proposed by \citet{Bessaih2015}, we considered the expected $\ell_2$-error norm, $\Lambda$, between 
downscaled and reference solutions, defined by

\begin{equation}
    \Lambda(t) = \mathbb{E} \left[ \int_\Omega \left( f(\bm{x}; t) - \hat{f}(\bm{x}; t) \right)^2 d \Omega \right], 
    \label{eqn:measure}
\end{equation}
where $\mathbb{E}$ is the average operator. Note that for CDA-downscaled ensembles, 
$\Lambda$ is proportional to the variance of the observational noise 
$\sigma^2_{obs}$ \citep{Bessaih2015}, and this theoretical estimate was computationally verified by \citet{HammoudCOMG2022}.
Unless specified otherwise and for the sake of brevity, the results presented in this study focus on the temperature RRMSE values, and the results corresponding to the velocity components and different $\mathcal{S}$, $\mathcal{T}$, and Ra number values are presented in the Supplementary Material.

\end{document}